\documentclass[a4paper]{article}

\usepackage{algorithm} 
\usepackage{algpseudocodex} 

\usepackage{amsmath} 
\usepackage{amssymb} 
\usepackage{amsthm} 
\usepackage[shortlabels]{enumitem}
    \setlist{nosep}
    \setlist[enumerate,1]{label={(\roman*)}}
\usepackage{todonotes}
\usepackage[font=footnotesize,margin=1cm]{caption}
\usepackage{subcaption}
\usepackage{graphicx}
    \graphicspath{{results/}}
\usepackage{mathtools}
\usepackage{orcidlink}

\usepackage{pgfplots}
\usepackage{tikz}
    \pgfplotsset{compat=newest}
    \usetikzlibrary{plotmarks}
    \usetikzlibrary{calc}
    \usetikzlibrary{positioning}
    \usetikzlibrary{arrows}
    \usetikzlibrary{arrows.meta}
    \usetikzlibrary{graphs}
    \usetikzlibrary{matrix}
    \usepgfplotslibrary{patchplots}
\usepackage{multicol}
\usepackage{xspace}
\usepackage{xcolor} 
\usepackage[all]{nowidow}
\usepackage{cite}
\usepackage[nameinlink]{cleveref} 

\hypersetup{
        colorlinks=false,
        pdfborder={0,0,0},
        pdfborderstyle={/S/U/W 0},
        }

\definecolor{darkgreen}{RGB}{15, 135, 0}  

\theoremstyle{plain}
\newtheorem{theorem}{Theorem}

\newtheorem{example}[theorem]{Example}
\newtheorem{lemma}[theorem]{Lemma}
\newtheorem{proposition}[theorem]{Proposition}
\newtheorem{remark}[theorem]{Remark}
\newtheorem{assumption}[theorem]{Assumption}
\theoremstyle{}

\makeatletter
\def\thm@space@setup{%
  \thm@preskip=8pt plus 2pt minus 2pt   
  \thm@postskip=6pt plus 2pt minus 2pt  
}
\makeatother

\newcommand{\norm}[1]{\Vert #1 \Vert}

\newcommand{\NN}{\mathbb{N}}
\newcommand{\RR}{\mathbb{R}}
\newcommand{\RRnonneg}{\RR_{\geq0}}

\newcommand{\Qbfc}{\hat{\mathsf{Q}}}
\newcommand{\Rbf}{\mathsf{R}}

\newcommand{\Qcost}{\mathbf{Q}}
\newcommand{\Rcost}{\mathbf{R}}

\newcommand{\und}{,~}
\newcommand{\eye}{I} 

\newcommand{\Uad}{\mathcal{U}_{\mathrm{ad}}}
\newcommand{\UadOutputFeedback}{\tilde{\mathcal{U}}_{\mathrm{ad}}}
\newcommand{\UadOptimalityFeedback}{\hat{\mathcal{U}}_{\mathrm{ad}}}
\newcommand{\inv}{^{-1}}
\newcommand{\tp}{^{\hspace{-1pt}\mathsf{T}}}

\AtBeginDocument{\renewcommand{\tilde}{\widetilde}}
\AtBeginDocument{\renewcommand{\hat}{\widehat}}

\DeclareMathOperator*{\rank}{rank}

\newcommand{\dt}{\,\mathrm{d}t}

\newcommand{\ds}{\,\mathrm{d}s}
\newcommand{\ddt}{\frac{\mathrm{d}}{\mathrm{d}t}}

\newcommand{\pH}{pH\xspace}
\newcommand{\portHamiltonian}{port-Ha\-mil\-to\-ni\-an\xspace}

\renewcommand{\vec}[1]{\begin{bmatrix} #1 \end{bmatrix}}
\newcommand{\textvec}[1]{[\begin{smallmatrix} #1 \end{smallmatrix}]}

\newcommand{\lebesgue}{L} 

\newcommand{\ham}{\mathcal{H}}
\newcommand{\hamc}{\hat{\mathcal{H}}} 
\newcommand{\etac}{\nabla \hamc} 
\newcommand{\Detac}{\nabla^2 \hamc} 
\renewcommand{\eta}{\nabla \ham}

\newcommand{\zc}{\hat{z}}
\newcommand{\dzc}{\dot{\hat{z}}}
\newcommand{\fc}{\hat{f}}
\newcommand{\gc}{\hat{g}}
\newcommand{\hc}{\hat{h}}
\newcommand{\Hc}{\hat{\mathsf{H}}}
\newcommand{\uc}{\hat{u}}
\newcommand{\yc}{\hat{y}}
\newcommand{\ucc}{\hat{\hat{u}}}
\newcommand{\ycc}{\hat{\hat{y}}}
\newcommand{\Pc}{P_{\mathrm{c}}}
\newcommand{\Pf}{P_{\mathrm{f}}}

\newcommand{\tPf}{\tilde{P}_{\mathrm{f}}}

\newcommand{\Jc}{\hat{J}}
\newcommand{\Rc}{\hat{R}}

\newcommand{\Mc}{\hat{M}}
\newcommand{\ellc}{\hat{\ell}}

\newcommand{\zluen}{\tilde{z}}
\newcommand{\dzluen}{\dot{\tilde{z}}}
\newcommand{\fluen}{\tilde{f}}
\newcommand{\gluen}{\tilde{g}}
\newcommand{\hluen}{\tilde{h}}
\newcommand{\uluen}{\tilde{u}}
\newcommand{\yluen}{\tilde{y}}
\newcommand{\nluen}{\tilde{n}}

\newcommand{\zekf}{\overline{z}}
\newcommand{\dzekf}{\dot{\overline{z}}}

\newcommand{\uekf}{\overline{u}}
\newcommand{\yekf}{\overline{y}}

\newcommand{\zk}{z_{1}}
\newcommand{\znk}{z_{2}}
\newcommand{\gk}{g_{1}}
\newcommand{\gnk}{g_{2}}
\newcommand{\rk}{r_{1}}
\newcommand{\rnk}{r_{2}}
\newcommand{\yk}{y_{1}}

\newcommand{\etabar}{\overline{\nabla}\ham}
\newcommand{\gbar}{\overline{g}}

\newcommand{\ubar}{\overline{u}}
\newcommand{\ybar}{\overline{y}}

\newcommand{\tolrel}{\varepsilon_{\text{rel}}} 
\newcommand{\tolabs}{\varepsilon_{\text{abs}}} 

\newcommand{\Omegatest}{\Omega_{\text{test}}}

\newcommand{\interconnectionmatrix}{\mathbf{F}}
\newcommand{\outputweighting}{\mathbf{W}}

\pgfplotsset{
  custom y log style/.style={
      yticklabel={
        \pgfkeys{/pgf/fpu=true}
        \pgfmathparse{exp(\tick)}%
        \pgfmathprintnumber[fixed,fixed zerofill, precision=2]{\pgfmathresult}
        \pgfkeys{/pgf/fpu=false}
      }
  }
}

\begin{document}

\title{
    Nonlinear systems and passivity: feedback control, model reduction, and time discretization%
    } 

\author{
    Tobias Breiten\,\orcidlink{0000-0002-9815-4897}\thanks{Institute of Mathematics, Technische Universität Berlin, Germany.}\,\,\thanks{\texttt{tobias.breiten@tu-berlin.de}}\,\,,
    Attila Karsai\,\orcidlink{0000-0001-9905-433X}\footnotemark[1]\,\,\thanks{\texttt{karsai@math.tu-berlin.de}}
}

\maketitle

\begin{abstract}                          
    Dynamical systems can be used to model a broad class of physical processes, and conservation laws give rise to system properties like passivity or \portHamiltonian structure. 
    An important problem in practical applications is to steer dynamical systems to prescribed target states, and feedback controllers combining a regulator and an observer are a powerful tool to do so. 
    However, controllers designed using classical methods do not necessarily obey energy principles, which makes it difficult to model the controller-plant interaction in a structured manner. 
    In this paper, we show that the combination of an optimal feedback law characterized by the Hamilton--Jacobi--Bellman equation and output feedback gives rise to passivity properties of the controller that are independent of the plant structure.
    Furthermore, we state conditions for the controller to have a \portHamiltonian realization and show that a model order reduction scheme can be deduced using the framework of nonlinear balanced truncation.
    To illustrate our results, we numerically realize the controller using the policy iteration
    and computationally verify passivity via a custom passivity-preserving discrete gradient scheme suitable for a wide class of passive systems.
\end{abstract}


\section{Introduction}

To enforce some behavior on the state of a dynamical system, a common approach is closed loop controller design, where measurements of the dynamical system to be controlled, the \emph{plant}, are fed into a second dynamical system, the \emph{controller}, which then feeds back to the plant.
In the context of modeling real-world phenomena, often additional structure is present in the system dynamics.
One example is passivity, which describes systems unable to generate energy on their own~\cite{willems72-dissipative1,willems72-dissipative2}.
If two passive systems are connected using \emph{power-conserving interconnection}, then the coupled system is again passive~\cite[Theorem 2.10]{sepulchre97-constructive}.
This fact allows for a modularized modeling approach of large-scale dynamical systems if all subsystems are passive~\cite{mehrmann23-control}.
Further, passivity properties can be incorporated in the algebraic structure of the dynamics by means of the \emph{\portHamiltonian} formalism~\cite{vanderschaft14-port}.

Unfortunately, already for linear systems, classical feedback control design methods like \emph{linear quadratic Gaussian} (LQG) control fail to construct a passive controller, even if the plant was passive to begin with, see, e.g.,~\cite{johnson79-state}.
As a remedy, in~\cite{breiten22-error} a modification of LQG control is discussed.
Here, the key idea is to modify the underlying cost functionals in such a way that passivity of the controller is enforced, and the controller turns out to be a combination of optimal state feedback and output feedback.
This line of research has not yet been explored for nonlinear systems, for which optimal feedback laws characterized by the \emph{Hamilton--Jacobi--Bellman equation}~\cite{crandall83-viscosity}.
Another prominent control scheme incorporating passivity is \emph{passivity based control}~\cite{ortega89-adaptive,ortega02-interconnection}, where typically \portHamiltonian systems are considered.
Since not all passive systems are \portHamiltonian, the latter can pose difficulties.

Since in practical applications, the underlying systems are often infinite-dimensional, discretization in the spatial variables typically leads to high-dimensional models.
Handling these models is computationally expensive.
Reduced order modeling offers methods to reduce the model complexity while maintaining accuracy.
Similar to the situation for controller design, classical approaches to construct reduced order models do not necessarily preserve passivity.
Structure-preserving methods for linear systems use tangential interpolation~\cite{gugercin12-structure}, parameter optimization~\cite{schwerdtner23-sobmor}, spectral factorization techniques~\cite{breiten22-passivity}, or system balancing~\cite{breiten22-error} based on modified cost functionals that are also used to construct controllers.
For nonlinear systems, apart from balancing techniques~\cite{kawano18-structure,sarkar23-structure}, also Petrov--Galerkin methods are suitable~\cite{buchfink23-symplectic,buchfink24-model,altmann25-structure}.

Our main contributions are listed below.
\begin{itemize}
\item
    We extend the controller results of~\cite{breiten22-error} to the nonlinear setting.
    In particular, combining an optimal feedback law with output feedback, we state a nonlinear controller that is passive independent of the plant structure.
    Unlike the procedure for linear systems, where the controller stems from modified cost functionals, we first state the controller and then discuss in what sense it is optimal.
    For that, we make use of classical results on \emph{inverse optimal control problems} as introduced by~\cite{kalman64-when} and studied in, e.g.,~\cite{moylan73-nonlinear} for nonlinear systems.
\item
    As in~\cite{breiten22-error} and motivated by~\cite{kawano18-structure,sarkar23-structure}, based on the value function underlying our passive controller design, we develop passivity-preserving nonlinear balancing techniques.
    In contrast to the approaches focused on \portHamiltonian systems, we consider general nonlinear passive systems represented in suitable coordinates, with the advantage that the representation is feasible under weaker assumptions on the system dynamics~\cite{karsai25-passivity}.
\item
    To numerically verify the passivity of our controller, we combine discrete gradients~\cite{harten83-upstream,itoh88-hamiltonian,gonzalez96-time} and the previously mentioned coordinate representation to state a passivity-preserving time discretization scheme.
    The method is similar to the one presented in~\cite{mclachlan99-geometric}.
\end{itemize}

The paper is structured as follows.
\Cref{sec:problem-setting} summarizes feedback control and passivity, states the optimal control problem underlying our controller design, and motivates the controller structure.
In \Cref{sec:passive-controller}, we present our passive controller, discuss when it steers the plant trajectory to the origin, and comment on the possibility of a \portHamiltonian realization.
In \Cref{sec:optimality-interpration}, we apply results from Hamilton--Jacobi--Bellman theory and inverse optimal control problems to interpret the controller structure as a modification of the associated cost functionals.
\Cref{sec:model-reduction} uses balanced truncation with two value functions to state a passivity-preserving model reduction scheme.

In \Cref{sec:numerical-experiments}, we showcase the effectiveness of the controller using numerical experiments. 
We study its behavior and verify the passivity of the controller using our passivity-preserving time discretization scheme based on discrete gradients.
Finally, in \Cref{sec:conclusion} we summarize our findings and give an outlook on possible topics for future research.

\paragraph*{Notation.}
We denote the set of all $k$-times continuously differentiable functions from $U$ to $V$ by $C^k(U,V)$ and abbreviate $C^k \coloneq C^k(U,V)$ when the spaces $U$ and $V$ are clear from context.
The Jacobian of a function~$f\colon \RR^n \to \RR^n$, and the gradient and Hessian of a function~$\ham\colon \RR^n \to \RR$ at a point $z$ are denoted by $Df(z)$, $\nabla \ham(z)$ and $\nabla^2 {\ham}(z)$, respectively.
Further, we write $A \succeq 0$ if $z\tp A z \geq 0$ for all $z\in \RR^{n}$, and $A \succ 0$ if $z\tp A z > 0$ for all $z\in\RR^{n}\setminus\{0\}$, where $\cdot\tp$ denotes the transpose.
The rank of the matrix~$A$ is denoted by $\rank(A)$.
Mostly, we suppress the time dependency of functions and write~$z$ instead of~$z(t)$.
For $p \in [1,\infty]$ we denote the associated Lebesgue space as $\lebesgue^p$.
For a given control function $u\colon [0,\infty) \to \RR^m$, we occasionally denote the associated trajectory~$z$ at time $t$ by $z(t;u)$.

\section{Problem setting}
\label{sec:problem-setting}

We consider a plant of the form 
\begin{equation}\label{eq:plant}
    \begin{aligned}
        \dot{z}(t) & = f(z(t)) + g(z(t)) u(t), ~~ z(0) = z_0,\\
        y(t) & = h(z(t))
    \end{aligned}
\end{equation}
with continuous functions $f\colon \RR^n \to \RR^n$, $g \colon \RR^n \to \RR^{n,m}$ and $h\colon \RR^n \to \RR^m$ satisfying $f(0) = 0$.
Here, $t \geq 0$ is the time, $z(t)\in\RR^n$ is the state of the system with initial condition $z(0) = z_0$, $u(t) \in \RR^m$ is a control input, and $y(t)\in\RR^m$ models measurements.
Throughout the paper, we assume existence of solutions of~\eqref{eq:plant} on $[0,\infty)$ in the sense of Carathéodory~\cite[Theorem 1.1]{coddington55-theory}.
Moreover, for the sake of brevity we mostly suppress the time dependency of~$z$,~$u$ and~$y$.

Our motivation is to find a control~$u$ that steers the associated trajectory~$z(t;u)$ of~\eqref{eq:plant} to the origin as~$t \to \infty$.
For this, we aim to use a second dynamical system, the \emph{controller} 
\begin{equation}\label{eq:abstract-controller}
    \begin{aligned}
        \dzluen & = \fluen(\zluen) + \gluen(\zluen) \uluen, ~~ \zluen(0) = \zluen_0,
        \\
        \yluen & = \hluen(\zluen),
    \end{aligned}
\end{equation}
connected to~\eqref{eq:plant} via the power-conserving interconnection as shown in \Cref{fig:feedback}, where $\interconnectionmatrix \in \RR^{m,m}$ is an invertible, symmetric interconnection matrix to be specified.
As before,~$\fluen\colon \RR^{\nluen} \to \RR^{\nluen}$,~$\gluen \colon \RR^{\nluen} \to \RR^{\nluen,m}$ and~$\hluen\colon \RR^{\nluen} \to \RR^m$ are assumed continuous, and~$\zluen_0 \in \RR^{\nluen}$.
Although in general, the state dimension of the controller $\nluen \in \NN$ may be different from the state dimension of the plant, in the following we focus on the case~$\nluen = n$.
The controller does not have access to the plant state~$z$, but only to measurements~$y$ of~$z$ multiplied by the interconnection matrix~$\interconnectionmatrix$.

\begin{figure}
\centering
\begin{tikzpicture}
    [
        every node/.style={
            rectangle, 
            rounded corners=2pt, 
            },
        thick,
        >=stealth,
    ]
    
    \node[fill=black!5!white] (plant) at (0,0.5) {plant~\eqref{eq:plant}};
    
    \node[fill=blue!8!white] (luenberger) at (0,-0.5) {controller~\eqref{eq:abstract-controller}};

    \draw[->] (plant.east) -- +(1,0) |- (luenberger.east);
    \draw[<-] (plant.west) -- +(-1,0) |- (luenberger.west);

    \node[fill=none,anchor=west] at ($(plant.east)+(1,-0.5)$) {$\uluen = \interconnectionmatrix y$};
    \node[fill=none,anchor=east] at ($(plant.west)+(-1,-0.5)$) {$u = - \interconnectionmatrix \yluen$};
\end{tikzpicture}
\caption{Plant~\eqref{eq:plant} and controller~\eqref{eq:abstract-controller} connected via power-conserving via the interconnection matrix~$\interconnectionmatrix$.}
\label{fig:feedback}
\end{figure}

Our goal is to ensure that the controller~\eqref{eq:abstract-controller} is \emph{passive}.
This is motivated by the fact that if both~\eqref{eq:plant} and~\eqref{eq:abstract-controller} are passive, then the interconnection of the two as depicted in \Cref{fig:feedback} remains passive, see, e.g.,~\cite[Proposition 2.11]{sepulchre97-constructive}.
We do, however, not make the passivity of~\eqref{eq:plant} a general assumption.

\subsection{Passivity}
For systems of the form~\eqref{eq:plant}, passivity takes the form of the existence of nonnegative \emph{storage function $\ham\colon \RR^n \to \RRnonneg$} satisfying the \emph{dissipation inequality}
\begin{equation}\label{eq:energy-inequality}
    \ham(z(t_1)) - \ham(z(t_0)) \leq \int_{t_0}^{t_1} y(t)\tp u(t) \dt
\end{equation}
for all $t_1 \geq t_0 \geq 0$ and controls~$u$ along solutions of~\eqref{eq:plant}.
Passivity can be characterized as follows.
If $\ham\colon\RR^n \to \RRnonneg$ is continuously differentiable and there exists $\ell\colon \RR^n \to \RR^p$ with 
\begin{equation}\label{eq:lure}
    \begin{aligned}
        \eta(z)\tp f(z) & = - \ell(z)\tp \ell(z) \\
        g(z)\tp \eta(z) & = h(z)
    \end{aligned}
\end{equation}
for all $z\in \RR^n$, then~\eqref{eq:plant} is passive and~$\ham$ is a storage function.
This follows from the \emph{energy balance}
\begin{equation}\label{eq:energybalance}
    \begin{aligned}
    &\ham(z(t_1)) - \ham(z(t_0)) 
    \\
    =& \int_{t_0}^{t_1} \ddt \ham(z) \dt 
    = \int_{t_0}^{t_1} \eta(z)\tp \dot{z} \dt 
    \\
    =& \int_{t_0}^{t_1} \eta(z)\tp f(z) + \eta(z)\tp g(z) u \dt 
    \\
    =& \int_{t_0}^{t_1} - \ell(z)\tp \ell(z) + y\tp u \dt.
    \end{aligned} 
\end{equation}
Under additional assumptions, the reverse implication is also true~\cite{moylan74-implications,hill80-dissipative}.

For modeling purposes it can be advantageous to have additional algebraic properties in the dynamics.
In this context the \portHamiltonian formalism is central~\cite{vanderschaft14-port,mehrmann23-control}.
Usually, the system~\eqref{eq:plant} with storage function~$\ham$ is called \emph{\portHamiltonian} if there exist functions $J, R \colon \RR^n \to \RR^{n,n}$ with $J(z) = - J(z)\tp$ and $R(z) = R(z)\tp \succeq 0$ and 
\begin{equation}\label{eq:ph-a}
    f(z) = (J(z) - R(z)) \eta(z)
\end{equation}
for all $z\in\RR^n$.
As discussed in, e.g.,~\cite{karsai25-passivity,ortega02-interconnection, wang03-generalized}, the problem of finding suitable operators $J$ and $R$ for a given storage function~$\ham$ is not trivial.
Another possibility to algebraically encode passivity are monotone structures, see, e.g.,~\cite{camlibel23-port,gernandt25-port}.
If we assume that $\eta$ is bijective, then~$r\coloneq-f\circ\eta\inv$ satisfies
\begin{equation}\label{eq:r-definition}
    f(z) = - r(\eta(z)) \und v\tp r(v) \geq 0
\end{equation}
for all $z,v \in \RR^n$.
Existence of functions~$r \colon \RR^n \to \RR^n$ satisfying~\eqref{eq:r-definition} may also be obtained using weaker assumptions, see~\cite[Proposition 6]{karsai25-passivity}.
Both representations~\eqref{eq:ph-a} and~\eqref{eq:r-definition} imply passivity since $\eta(z)\tp f(z) \leq 0$ for all $z \in \RR^n$ due to the respective properties of~$R$ and~$r$.


\subsection{The optimal control problem}

If access to the state~$z$ of~\eqref{eq:plant} is available via $y=z$, then one possibility to construct a controller is to consider an optimal control problem and choose the control~$u$ as the associated optimal feedback law.
We will focus on the problem
\begin{equation}\label{eq:ocp}
    \hamc(z_0) = \inf_{u \in \Uad(z_0)} \left\{ \frac{1}{2} \int_{0}^{\infty} y\tp \Qcost y + u\tp \Rcost u \dt \right\}
\end{equation}
subject to the dynamics~\eqref{eq:plant} with the initial condition $z(0) = z_0$ and associated value function $\hamc \colon \RR^n \to \RRnonneg$.
Here, $\Qcost, \Rcost \in \RR^{n,n}$ are symmetric positive semidefinite weighting matrices with~$\Rcost$ invertible, and the set of admissible controls $\Uad(z_0)$ is defined as
\begin{equation*}
    \Uad(z_0) \coloneq
    \{ u \in \lebesgue^{2}([0,\infty), \RR^m) ~|~ \text{$u$ renders the integral in~\eqref{eq:ocp} finite}\}. 
\end{equation*}
We refer the reader to, e.g.,~\cite{fleming06-controlled} for a detailed discussion on the existence of optimal controls for~\eqref{eq:ocp}, and make the following simplifying assumption throughout the paper.

\begin{assumption}\label{as:smoothness-hamc}
    The set of admissible controls $\Uad(z_0)$ is nonempty for all initial values $z_0\in\RR^n$, and the value function $\hamc$ of~\eqref{eq:ocp} is continuously differentiable.
    Moreover, the feedback law $u^*(z) = - \Rcost\inv g(z)\tp \etac(z)$ is an optimal control for~\eqref{eq:ocp}.
\end{assumption}


The optimality of the feedback law $u^*(z) = - \Rcost\inv g(z)\tp \etac(z)$ is a classical result from Hamilton--Jacobi--Bellman theory~\cite{crandall83-viscosity}.
For a formal derivation, we start with \emph{Bellman's dynamic programming principle}, which follows from the semigroup property of~\eqref{eq:plant} and states that for any $T>0$ we have
\begin{equation*}
    0 
    = 
    \inf_{u \in \Uad(z_0)}
    \Bigg\{
    \frac12 \int_{0}^{T} y\tp \Qcost y + u\tp \Rcost u \dt 
    + 
    \hamc(z(T; u))
    - 
    \hamc(z_0)
    \Bigg\}
    .
\end{equation*}
Dividing by~$T$, we obtain 
\begin{equation*}
    0
    = 
    \inf_{u \in \Uad(z_0)} 
    \Bigg\{
    \frac{1}{2T}
    \int_{0}^{T} 
        y\tp \Qcost y + u\tp \Rcost u
    \dt
    + \frac{\hamc(z(T; u)) - \hamc(z_0)}{T} 
    \Bigg\}.
\end{equation*}
Under smoothness assumptions, letting $T \to 0$ and using $\lim_{T \to 0} \tfrac1T \int_{0}^{T} \phi(t) \dt = \phi(0)$ for continuous~$\phi$ we arrive at the \emph{dynamic programming equation}%
\begin{equation}\label{eq:control-dpe}
    \inf_{u \in \RR^m}
    \big\{
    \etac(z)\tp (f(z) + g(z) u)
    +
    \tfrac12
    h(z)\tp \Qcost h(z) 
    + 
    \tfrac12
    u\tp \Rcost u
    \big\}
    =
    0
\end{equation}
for all $z \in \RR^n$, where for notational convenience we have replaced~$z_0$ by~$z$, $y(0)$ by $h(z)$ and~$u(0)$ by~$u$.
Since the term is quadratic in~$u$, we can explicitly compute the minimizer and obtain the feedback law~$u^*(z) = - \Rcost\inv g(z)\tp \etac(z)$.

Unfortunately, the state $z$ is not always available in applications, so that directly implementing~$u^*$ may not be feasible.
In this case, one possibility is to resort to an approximation~$\zluen$ of~$z$ and choosing $u^*(\zluen) = - \Rcost\inv g(\zluen)\tp \etac(\zluen)$ as the control.
A common approach is choosing~$\zluen$ as the state of a Luenberger observer~\cite{luenberger64-observing}, which copies the original dynamics and adds an innovation term correcting for deviances using observations.
Plugging~$u^*(\zluen)$ in~\eqref{eq:plant} and adding an innovation term yields the controller
\begin{equation}\label{eq:luenberger}
    \begin{aligned}
        \dzluen & = f(\zluen) - g(\zluen) \Rcost\inv g(\zluen)\tp \etac(\zluen)
        + L(\cdot, \zluen)(\interconnectionmatrix y - \outputweighting h(\zluen)), 
        \\
        \yluen & = \interconnectionmatrix\inv \Rcost\inv g(\zluen)\tp \etac(\zluen),
    \end{aligned}
\end{equation}
where $\interconnectionmatrix y - \outputweighting h(\zluen)$ is the innovation term with invertible weighting matrix $\outputweighting \in \RR^{m,m}$, $L(t, \zluen(t)) \in \RR^{n,m}$ is the \emph{observer gain} at time~$t$, and $\zluen(0) = \zluen_0$ is an initial condition.
If we interpret $\interconnectionmatrix y = \uluen$ as a control input to~\eqref{eq:luenberger}, then the feedback law $u(\zluen) = - \Rcost\inv g(\zluen)\tp \etac(\zluen) = - \interconnectionmatrix \yluen$ is the result of the power-conserving interconnection of~\eqref{eq:plant} and~\eqref{eq:luenberger} via the interconnection matrix $\interconnectionmatrix$, see \Cref{fig:feedback}.


The choice~$\outputweighting = \interconnectionmatrix$ is particularly interesting, since then the innovation term directly corresponds to the difference between observation and estimated output.

\section{The passive controller}
\label{sec:passive-controller}

In this section, we show that the choices $\interconnectionmatrix = \Rcost^{-1/2}$, $\outputweighting = \Qcost^{1/2}$ and $L(\cdot,\zluen) = g(\zluen) \Rcost^{-1/2}$ in~\eqref{eq:luenberger} leads to a passive controller for the plant~\eqref{eq:plant}.
Furthermore, we remark on \portHamiltonian formulations of the controller.

\begin{theorem}[passive controller]\label{thm:passive-controller}
    Let \Cref{as:smoothness-hamc} hold.
    Then the system
    \begin{equation}\label{eq:passive-controller}
        \begin{aligned}
        \dzc 
        & 
        = 
        f(\zc) - g(\zc) \Rcost\inv g(\zc)\tp \etac(\zc) 
        +  g(\zc) \Rcost^{-1/2} (\uc - \Qcost^{1/2} h(\zc)),
        \\
        \yc & = \Rcost^{-1/2} g(\zc)\tp \etac(\zc)
        \end{aligned}
    \end{equation}
    with the initial condition $\zc(0) = \zc_0$ is passive with storage function $\hamc \geq 0$.
\end{theorem}
\begin{proof}
    Plugging $u^*(z) = - \Rcost\inv g(z)\tp \etac(z)$ in~\eqref{eq:control-dpe}, we obtain the \emph{Hamilton--Jacobi--Bellman equation}
    \begin{equation}\label{eq:control-hjb}
        \etac(z)\tp f(z) - \tfrac12 \etac(z)\tp g(z) \Rcost\inv g(z)\tp \etac(z)
        + \tfrac12 h(z)\tp \Qcost h(z) = 0.
    \end{equation}
    Since the running cost in~\eqref{eq:ocp} is nonnegative, it immediately follows that $\hamc \geq 0$.
    Let us denote $\gc \coloneq g \Rcost^{-1/2}$ and $\fc \coloneq f - \gc \gc\tp \etac - \gc \Qcost^{1/2} h$ such that the dynamics of~\eqref{eq:passive-controller} read as $\dzc = \fc(\zc) + \gc(\zc) \uc$.
    To verify the passivity of~\eqref{eq:passive-controller}, we show that~\eqref{eq:lure} holds with~$f$,~$\eta$ and~$y$ replaced by~$\fc$,~$\etac$ and~$\yc = \gc\tp \etac$, respectively.
    For the first equation of~\eqref{eq:lure}, we notice that due to~\eqref{eq:control-hjb} we have
    \begin{align*}
        & \hphantom{=}~
        \etac(\zc)\tp \fc(\zc) 
        \\
        & 
        = 
        \etac(\zc)\tp 
        \left( 
            f(\zc) 
            - \gc(\zc) \gc(\zc)\tp \etac(\zc) 
            - \gc(\zc) \Qcost^{1/2} h(\zc) 
        \right) 
        \\
        & 
        = 
        -\tfrac12 h(\zc)\tp \Qcost h(\zc) 
        + \tfrac12 \etac(\zc)\tp g(\zc) \Rcost\inv g(\zc)\tp \etac(\zc) 
        \\
        & 
        ~~~\,
        - \etac(\zc)\tp \gc(\zc) \gc(\zc)\tp \etac(\zc)  
        - \etac(\zc)\tp \gc(\zc) \Qcost^{1/2} h(\zc) 
        \\
        & 
        = 
        -\tfrac12 h(\zc)\tp \Qcost h(\zc) 
        - \tfrac12 \etac(\zc)\tp \gc(\zc) \gc(\zc)\tp \etac(\zc) 
        \\
        & 
        ~~~\, 
        - \tfrac12 
        \left( 
            \etac(\zc)\tp \gc(\zc) \Qcost^{1/2} h(\zc) 
            + h(\zc)\tp \Qcost^{1/2} \gc(\zc)\tp \etac(\zc) 
        \right) 
        \\
        & 
        = 
        -\tfrac12 
            (
            \Qcost^{1/2} h(\zc) 
            + \gc(\zc)\tp \etac(\zc)
            )\tp
            (
            \Qcost^{1/2} h(\zc) 
            + \gc(\zc)\tp \etac(\zc)
            )
    \end{align*}
    so that $\ellc(\zc) = \frac{1}{\sqrt{2}} (\Qcost^{1/2}h(\zc) + \gc(\zc)\tp \etac(\zc))$ satisfies $\etac\tp \fc = - \ellc\tp \ellc$.
    Finally, we note that the second equation in~\eqref{eq:lure} is satisfied by construction so that passivity of~\eqref{eq:passive-controller} follows.
\end{proof}


As mentioned before, to recover a more commonly used innovation term, we need to consider~$\Qcost^{1/2} = \outputweighting = \interconnectionmatrix = \Rcost^{-1/2}$.
This coupling of $\Qcost$ and $\Rcost$ eliminates parts of the freedom in the choice of the weighting matrices.

At this point, it is not clear whether the interconnection of~\eqref{eq:plant} and~\eqref{eq:passive-controller} steers~$z(t) \to 0$ as~$t \to \infty$.
Since in the literature often $\interconnectionmatrix = \eye$ is considered to study interconnections, we note that the controller~\eqref{eq:passive-controller} can be rewritten as 
\begin{equation}\label{eq:passive-controller-rewritten}
    \begin{aligned}
    \dzc 
    & 
    = 
    f(\zc) 
    - g(\zc) \Rcost\inv g(\zc)\tp \etac(\zc)
    - g(\zc) \Rcost^{-1/2} \Qcost^{1/2} h(\zc)
    +  g(\zc) \Rcost^{-1} \ucc,
    \\
    \ycc & = \Rcost^{-1} g(\zc)\tp \etac(\zc),
    \end{aligned}
\end{equation}
where $\ucc = y$ and $u = -\ycc$ are the controller inputs and outputs for the interconnection matrix~$\interconnectionmatrix = \eye$, respectively.
To study the closed loop performance of~\eqref{eq:plant} and~\eqref{eq:passive-controller-rewritten} standard results such as small gain theorems may be applied, see, e.g.,~\cite[Section 2.1]{vanderschaft17-l2gain}.
Since we are particularly interested in the interconnection of passive systems, we adapt~\cite[Theorem 10.8.1]{isidori99-nonlinear} as follows.
Below, we write~\eqref{eq:plant}$\circ$\eqref{eq:passive-controller-rewritten} to denote the open loop interaction of~\eqref{eq:plant} and \eqref{eq:passive-controller-rewritten}, that is, taking $\ucc = y$ only.

\begin{proposition}[{\!\!\cite[Theorem 10.8.1]{isidori99-nonlinear}}]\label{prop:interconnection}
    In addition to \Cref{as:smoothness-hamc}, let~\eqref{eq:plant} be passive with storage function~$\ham$.
    Moreover, let symmetric matrices $\Qbfc, \Rbf \in \RR^{m,m}$ be given such that $\Qbfc + \Rbf$ is negative definite and
    \begin{align*}
        \ham(z(t_1)) - \ham(z(t_0)) & \leq \int_{t_0}^{t_1} y\tp u + u\tp \Rbf u \dt,
        \\
        \hamc(\zc(t_1)) - \hamc(\zc(t_0)) & \leq \int_{t_0}^{t_1} \ycc\tp \Qbfc \ycc + \ycc\tp \ucc \dt 
    \end{align*}
    for all $t_1 \geq t_0 \geq 0$ and controls~$u$ and~$\ucc$ along solutions of~\eqref{eq:plant} and~\eqref{eq:passive-controller-rewritten}, respectively. 
    Assume further that~\eqref{eq:plant}$\circ$\eqref{eq:passive-controller-rewritten} is zero-state detectable.
    Then the power-conserving interconnection of~\eqref{eq:plant} and~\eqref{eq:passive-controller-rewritten} is globally asymptotically stable.
\end{proposition}

\begin{proof}
    Consider $a = 1$, $S_1 = S_2 = \tfrac12 \eye$, $Q_1 = 0$, $R_1 = \Rbf$, $Q_1 = \Qbfc$ and $R_2 = 0$ in~\cite[Theorem 10.8.1]{isidori99-nonlinear}
\end{proof}

The following consequence of~\cite[Proposition 10]{karsai25-passivity}, gives sufficient conditions for the controller~\eqref{eq:passive-controller} to have a \portHamiltonian representation in the form~\eqref{eq:ph-a}.
As with the passivity of~\eqref{eq:passive-controller}, the result is independent of the structure of~\eqref{eq:plant}.
Here and in the following, we set $\fc \coloneq f - g g\tp \etac - g h$.

\begin{proposition}\label{thm:ph-representation-passive-controller}
    Assume $\fc \in C^1$, $\hamc \in C^2$ and let $\etac$ be bijective with $(\fc \circ \etac\inv)(0) = 0$.
    Further, let $\Detac(\zc)$ be invertible with $D\fc(\zc)  \Detac(\zc)\inv \preceq 0$ for all $\zc\in\RR^n$.
    Then~\eqref{eq:passive-controller} is \portHamiltonian with  $\fc(\zc) = (\Jc(\zc) - \Rc(\zc))\etac(\zc)$, where $\Jc(\zc) \coloneq \tfrac12(\Mc(\zc)- \Mc(\zc)\tp)$ and $\Rc(\zc) \coloneq -\tfrac12 (\Mc(\zc) + \Mc(\zc)\tp)$ with $\Mc(\zc) \coloneq \int_{0}^{1} D\fc(s\zc) \Detac(s\zc)\inv \ds$.
\end{proposition}

For linear systems $D\fc(\zc) \Detac(\zc)\inv \preceq 0$ follows from the equivalence of the Hamilton--Jacobi--Bellman equation~\eqref{eq:control-hjb} and a corresponding algebraic Riccati equation, see~\cite[Theorem 1]{breiten22-error}.
In the general nonlinear case, the assumption is not easily verified.

\section{Optimality interpretation}
\label{sec:optimality-interpration}


In the linear case discussed in~\cite{breiten22-error}, the controller structure~\eqref{eq:passive-controller} stems from choosing specific cost functionals in the LQG control and observation problems.
The goal in this section is to obtain similar results for the nonlinear case.
For the sake of brevity, we focus on the case~$\Qcost = \Rcost = \eye$, but note that the results can be adapted to more general weighting matrices.
The uncontrolled dynamics of~\eqref{eq:passive-controller} read 
\begin{equation*}
    \dzc = f(\zc) - g(\zc)g(\zc)\tp \etac(\zc) - g(\zc) h(\zc).
\end{equation*}
We aim to give both the output feedback~$-h$ and the combined feedback $- g\tp \etac - h$ interpretations in an optimality sense.
For this, we focus on the case that~\eqref{eq:plant} is passive in the sense of~\eqref{eq:lure}.
This can be thought of as an \emph{inverse optimal control problem} as introduced by~\cite{kalman64-when} and studied and applied to robust control design in, e.g.,~\cite{moylan73-nonlinear,sepulchre97-constructive,freeman96-robust}.
Inverse optimality concepts have been combined with observers in, e.g.,~\cite{jin21-inverse,zanon21-constrained} for the discrete time setting.


\subsection{Output feedback}

The next result states that in this case the output feedback term $-h = - g\tp \eta$ in~\eqref{eq:passive-controller} is the optimal control for a cost functional associated with the passivity structure of~\eqref{eq:plant}.
The result is a result on inverse optimal control problems and here formulated as a consequence of~\cite[Theorem 8.3]{haddad08-nonlinear}. 
Alternatively, the result may be obtained using verification theorems, see, e.g.,~\cite[Theorem I.7.1]{fleming06-controlled}.
Here and in the following, as in Lyapunov stability theory~\cite[Section 3.2]{haddad08-nonlinear}, a function~$\phi \colon \RR^n \to \RR$ is called \emph{radially unbounded} if $\norm{z} \to \infty$ implies $\phi(z) \to \infty$.

\begin{theorem}[interpretation of output feedback]\label{thm:interpretation-output-feedback}
    Assume that \eqref{eq:plant} is passive in the sense of~\eqref{eq:lure} with~$\ell \in C^0$ and~$\ham \in C^1$ being radially unbounded such that $\ham(0) = 0$, $\ham(z) > 0$ for $z\neq0$.
    Additionally, let
    \begin{equation}\label{eq:haddad-condition}
        \ell(z)\tp \ell(z) + h(z)\tp h(z) > 0
    \end{equation}
    for $z \neq 0$.
    Then $\ham$ is the value function of the optimal control problem 
    \begin{equation}\label{eq:interpretation-output-feedback-problem}
        \inf_{u \in \UadOutputFeedback(z_0)} 
        \left\{ 
            \frac12 \int_{0}^{\infty} 
            h(z)\tp h(z)
            + 2 \ell(z)\tp \ell(z) 
            + u\tp u \dt 
        \right\}
    \end{equation}
    subject to the dynamics~\eqref{eq:plant}, where
    \begin{multline*}
        \UadOutputFeedback(z_0) 
        \coloneq
        \{ u \in \lebesgue^{2}([0,\infty), \RR^m) ~|~ \text{$u$ renders the integral in~\eqref{eq:interpretation-output-feedback-problem} finite}
        \\
        \text{and} \lim_{t\to\infty} z(t; u) = 0 \},
    \end{multline*}
    and the feedback law $\tilde{u}^*(z) = - g(z)\tp \eta(z) = - h(z)$ is an optimal control.
\end{theorem}
\begin{proof}
    The result follows from~\cite[Theorem 8.3]{haddad08-nonlinear}.
    More precisely, consider $G \coloneq g$, $V \coloneq \ham$, $L_1 \coloneq h\tp h + 2 \ell\tp \ell$, $L_2 \coloneq 0$,  and $R_2 \coloneq \tfrac12\eye$.
    Then, assumptions~(8.59)--(8.61) and~(8.63) from~\cite[Theorem 8.3]{haddad08-nonlinear} are satisfied.
    For~(8.62), note 
    that~\eqref{eq:lure} together with~\eqref{eq:haddad-condition} implies 
    \begin{equation*}
        \eta(z)\tp f(z) - \eta(z)\tp g(z) g(z)\tp \eta(z) < 0
    \end{equation*}
    for $z \neq 0$.
    Then, the feedback law~$u^*$ is the feedback law~$\phi$ in~\cite[Theorem 8.3]{haddad08-nonlinear}.
\end{proof}

The above result is connected to the results of~\cite{moylan73-nonlinear} as follows.
Let~\eqref{eq:plant} be passive in the sense of~\eqref{eq:lure} with~$\ham(0) = 0$, and let Assumptions~1 to~3 from~\cite{moylan73-nonlinear} be satisfied for~\eqref{eq:plant} and~$m = h\tp h + 2 \ell\tp \ell$.
In~\cite[Theorem 3]{moylan73-nonlinear}, it is shown that the feedback law $-h = - g\tp \eta$ is optimal for the cost functional~$\ham$ if and only if the \emph{return difference condition} 
\begin{equation}\label{eq:rdc}
    \int_{0}^{\infty} (u + h(z))\tp (u + h(z)) \dt
    \geq
    \int_{0}^{\infty} u\tp u \dt 
\end{equation}
is satisfied for all asymptotically stabilizing $u \in L^2([0,\infty), \RR^m)$ subject to the dynamics~\eqref{eq:plant} with initial condition $z_0 = 0$.
For such~$u$, with $\ham(z(\infty)) = \ham(z_0) = 0$ and the energy balance~\eqref{eq:energybalance} we have 
\begin{equation*}
    \int_{0}^{\infty} h(z)\tp u \dt 
    =
    \int_{0}^{\infty} \ell(z)\tp \ell(z) \dt,
\end{equation*}
so that~\eqref{eq:rdc} is equivalent to 
\begin{equation*}
    \int_{0}^{\infty} 2 \ell(z)\tp \ell(z) + h(z)\tp h(z) \dt \geq 0,
\end{equation*}
which is clearly satisfied.

As the next remark illustrates, condition~\eqref{eq:haddad-condition} can be restrictive. 

\begin{remark}\label{rem:haddad-condition}
    If~\eqref{eq:plant} is a linear time-invariant passive system with quadratic storage function~$\ham(z) = \tfrac12 z\tp Q z$, $Q = Q\tp \succ 0 \in \RR^{n,n}$, then it can be written as a linear time-invariant \portHamiltonian system of the form 
    \begin{equation}\label{eq:lti-ph}
        \begin{aligned}
            \dot{z} & = (J-R)Qz + B u, \\
            y & = B\tp Q z,
        \end{aligned}
    \end{equation}
    where $J,R \in \RR^{n,n}$ satisfy $J = -J\tp$ and $R = R\tp \succeq 0$, see, e.g.,~\cite[Theorem 3]{willems72-dissipative2}.
    Then $\ell(z)\tp \ell(z) = - \eta(z)\tp f(z) = z\tp QRQ z$ so that~\eqref{eq:haddad-condition} reads as 
    \begin{equation*}
        z\tp Q R Q z + z\tp Q B B\tp Q z > 0
    \end{equation*}
    for $z \neq 0$.
    Since $Q$ is assumed invertible, this is equivalent to $R + BB\tp \succ 0$ or alternatively~$\rank([R \ B]) = n$.
    As discussed in, e.g.,~\cite[Proposition 5]{breiten25-structure}, this is sufficient but not necessary for the stabilizability and detectability of~\eqref{eq:lti-ph}.
\end{remark}

The next remark highlights the connection of our results to the ones obtained for the linear case~\cite{breiten22-error}.

\begin{remark}\label{rem:interpretation-output-feedback-linear}
    As in \Cref{rem:haddad-condition}, assume that~\eqref{eq:plant} is in the form~\eqref{eq:lti-ph} and denote $A \coloneq (J-R)Q$ and $C \coloneq B\tp Q$.
    In addition, assume that the system is stabilizable and detectable.
    Then the value function $\hamc$ of~\eqref{eq:ocp} is given by $\hamc(z_0) = \tfrac12 z_0\tp \Pc z_0$, where $\Pc = \Pc\tp \succ 0$ is the stabilizing solution to the algebraic Riccati equation 
    \begin{equation}\label{eq:control-are}
        A\tp \Pc + \Pc A - \Pc B B\tp \Pc + C\tp C = 0.
    \end{equation}
    We then have $\etac(\zc) = \Pc \zc$ and the controller~\eqref{eq:passive-controller} is given by 
    \begin{align*}
        \dzc & = (A - B B\tp \Pc - B C)\zc + B \uc, \\
        \yc & = B\tp \Pc \zc.
    \end{align*}
    Further, plugging the feedback $\tilde{u}^*(z) = -B\tp Q z$ into the dynamic programming equation associated with~\eqref{eq:interpretation-output-feedback-problem}, we obtain 
    \begin{equation*}
        z\tp Q A z - \tfrac12 C\tp C z + \tfrac12 z\tp Q B B\tp Q z + z\tp QRQ z  = 0
    \end{equation*}
    for all $z \in \RR^n$, which implies 
    \begin{equation*}
        Q A + A\tp Q - C\tp C + Q B B\tp Q + 2 Q R Q = 0.
    \end{equation*}
    Multiplying the above equation by $\tPf \coloneq Q\inv$ from the left and right, we obtain
    \begin{equation*}
        A \tPf + \tPf A\tp - \tPf C\tp C \tPf + B B\tp + 2R = 0.
    \end{equation*}
    This algebraic Riccati equation is also found in~\cite{breiten22-error} and can be understood as a modification of the filter equation from LQG control, which reads as 
    \begin{equation*}
        A \Pf + \Pf A\tp - \Pf C\tp C \Pf + B B\tp = 0.
    \end{equation*}
    Note also that~\eqref{eq:interpretation-output-feedback-problem} has similarities to the optimal control problem mentioned in~\cite[Remark 3]{breiten22-error}.
\end{remark}

\subsection{Combined feedback law}

If~\eqref{eq:plant} is passive in the sense of~\eqref{eq:lure}, then with $\uc = 0$ in~\eqref{eq:passive-controller} the dynamics of~\eqref{eq:passive-controller} become $\dot{\zc} = f(\zc) - g(\zc) g(\zc)\tp (\etac(\zc) + \eta(\zc))$.
To obtain an optimality result for the combined feedback law $-g\tp (\etac + \eta)$, we need to make the additional assumption that also $\ham + \hamc$, which is a storage function for the closed loop system, can serve as a storage function for~\eqref{eq:plant}.
To avoid confusion, we distinguish~$\hc \coloneq g\tp \etac$ from~$\yc$.

\begin{theorem}[optimality of feedback law]\label{thm:interpretation-optimality}
    Let \Cref{as:smoothness-hamc} hold and assume
    that~\eqref{eq:plant} is passive with storage function~$\ham \in C^1$.
    Further, assume that $(\etac(z) + \eta(z))\tp f(z) = - \ellc(z)\tp \ellc(z)$ for some $\ellc \in C^0$, and that $\ham + \hamc$ is radially unbounded with $(\ham + \hamc)(0) = 0$ and $(\ham + \hamc)(z) > 0$ for $z \neq 0$.
    If 
    \begin{equation}\label{eq:haddad-condition-2}
        \ellc(z)\tp \ellc(z) + h(z)\tp h(z) + \hc(z)\tp \hc(z) > 0
    \end{equation}
    for $z \neq 0$, then $\ham + \hamc$ is the value function of the optimal control problem
    \begin{equation}\label{eq:interpretation-optimality-problem}
            \inf_{u \in \UadOptimalityFeedback(z_0)} 
        \bigg\{ 
            \frac12 \int_{0}^{\infty} 
            h(z)\tp h(z) + \hc(z)\tp \hc(z)
            + 2 \ellc(z)\tp \ellc(z) 
            + u\tp u 
            \dt 
        \bigg\}
    \end{equation}
    subject to the dynamics~\eqref{eq:plant}, where 
    \begin{multline*}
        \UadOptimalityFeedback(z_0) \coloneq
        \{ u \in \lebesgue^{2}([0,\infty), \RR^m) ~|~ \text{$u$ renders the integral in~\eqref{eq:interpretation-optimality-problem} finite}  
        \\ \text{and} \lim_{t\to\infty} z(t; u) = 0 \},
    \end{multline*}
    and the feedback law $\hat{u}^*(z) = - g(z)\tp (\eta(z) + \etac(z)) = - (h(z) + \hc(z))$ is an optimal control.
\end{theorem}
\begin{proof}
    As before, the result is a consequence of~\cite[Theorem 8.3]{haddad08-nonlinear}.
    Consider $G \coloneq g$, $V \coloneq \ham + \hamc$, $L_1 \coloneq h\tp h + \hc\tp \hc + 2 \ell\tp \ell$, $L_2 \coloneq 0$,  and $R_2 \coloneq \tfrac12\eye$.
    Then, assumptions~(8.59)--(8.61) and~(8.63) from~\cite[Theorem 8.3]{haddad08-nonlinear} are satisfied.
    Assumption~(8.62) follows from~\eqref{eq:haddad-condition-2}.
\end{proof}

As with \Cref{thm:interpretation-output-feedback}, also \Cref{thm:interpretation-optimality} can be connected to~\cite{moylan73-nonlinear}.
For the sake of brevity, we omit the details.

As the next remark shows, the existence of $\hat{\ell}$ as in \Cref{thm:interpretation-optimality} can be a restrictive assumption.

\begin{remark}\label{rem:optimality-interpretation-linear}
    For linear time invariant \pH systems as in \Cref{rem:interpretation-output-feedback-linear}, the condition $(\etac(z) + \eta(z))\tp f(z) = - \hat{\ell}(z)\tp \hat{\ell}(z)$ in \Cref{thm:interpretation-optimality} translates to 
    \begin{equation}\label{eq:optimality-interpretation-condition-linear}
        A\tp (\Pc + Q) + (\Pc + Q) A \preceq 0,
    \end{equation}
    where $A = (J-R)Q$ for $J=-J\tp$, $R=R\tp \succeq 0$ and $Q=Q\tp \succ 0$, and $\Pc$ is the solution to~\eqref{eq:control-are} with $C = B\tp Q$.
    The inequality~\eqref{eq:optimality-interpretation-condition-linear} does not hold in general.
    To see this, consider the \portHamiltonian system defined by
    \begin{equation*}
        J = \vec{0 & -1 \\ 1 & 0},~
        R = \vec{1 & 0 \\ 0 & 0},~
        Q = \vec{1 & 0 \\ 0 & 1},~
        B = \vec{1 \\ -1}.
    \end{equation*}
    Then the stabilizing solution $\Pc$ of~\eqref{eq:control-are} with $C = B\tp Q$ is $\Pc = \textvec{\sqrt{2}-1 & 0 \\ 0 & 1}$.
    The eigenvalues of $A\tp (\Pc + Q) + (\Pc + Q) A$ are $\lambda_{1,2}= -\sqrt{2} \pm 2 \sqrt{2 - \sqrt{2}}$ with $\lambda_1 > 0$ and $\lambda_2 < 0$.
\end{remark}

\section{Application to model order reduction}
\label{sec:model-reduction}

In practical applications, the simulation of the full order dynamics~\eqref{eq:plant} can be prohibitively expensive.
In this regard, model order reduction schemes can be useful.
A common approach for linear time-invariant systems is \emph{balanced truncation}~\cite{antoulas05-approximation}, where system states that are not of interest are discarded.
This is done by first transforming the state space to appropriate coordinates (\emph{balancing}) and then truncating the transformed state space vector (\emph{truncation}).
The idea has been extended to nonlinear systems in~\cite{scherpen93-balancing}.
As mentioned in the introduction, structural properties like passivity are generally lost in this procedure.
For linear systems, balancing with respect to Gramians also used to construct passive controllers leads to a passivity-preserving model reduction method~\cite{breiten22-error}.
For nonlinear \portHamiltonian systems, similar ideas have been used in~\cite{kawano18-structure,sarkar23-structure}.
Our goal in this section is to present results like~\cite{breiten22-error} for nonlinear passive systems that are not necessarily \portHamiltonian.
We use techniques similar to~\cite{kawano18-structure,sarkar23-structure} and show that if~\eqref{eq:plant} is passive, then under smoothness assumptions we show that there exists a coordinate transformation that renders the storage function~$\ham$ of~\eqref{eq:plant} quadratic and puts the value function~$\hamc$ in a simple algebraic form.
In contrast to~\cite{kawano18-structure,sarkar23-structure}, we do not assume the dynamics to be $C^\infty$ and state the required smoothness explicitly.

We start with the following consequence of Morse theory, which appeared similarly in~\cite[Lemmas 4.1 and 4.2]{scherpen93-balancing}.

\begin{lemma}\label{lem:morse}
    Let~$p, q \colon \RR^n \to \RR$ be $C^{k+2}$ and $C^{k}$ functions with $k \geq 2$, respectively, that satisfy $p(0) = q(0) = 0$, $\nabla p(0) = 0$ and $\nabla^2 p(0), \nabla^2 q(0) \succ 0$.
    Then, there exists an open set
    $U \subseteq \RR^n$ with~$0 \in U$ and a $C^k$-isomorphism $\phi \colon U \to \phi(U),~ z \mapsto \bar{z}$ 
    with $\phi(0) = 0$
    such that on $\phi(U)$ we have
    \begin{gather*}
        \bar{p}(\bar{z}) \coloneq (p \circ \phi\inv)(\bar{z}) = \frac12 \bar{z}\tp \bar{z},
        \\
        \bar{q}(\bar{z}) \coloneq (q\circ\phi\inv)(\bar{z}) = \frac12 \bar{z}\tp \mathsf{Q}(\bar{z}) \bar{z}\tp
    \end{gather*}
    for a pointwise symmetric positive definite $C^{k-2}$ function $\mathsf{Q} \colon \phi(U) \to \RR^{n,n}$ with $\mathsf{Q}(0) = \nabla^2 q(0) \succ 0$.
\end{lemma}
\vspace{-15pt}
\begin{proof}
    The existence of an open set~$U$ and a $C^{k}$-isomorphism $\phi \colon U \to \RR^n, ~ z \mapsto \bar{z}$ with~$\phi(0)=0$ and~$\bar{p}(\bar{z}) = (p \circ \phi^{-1})(\bar{z}) = \tfrac12 \bar{z}\tp \bar{z}$ follows from Morse's lemma~\cite[Theorem 4 in Section 7]{lang85-differential}.
    Applying~\cite[Lemma 2.1]{milnor16-morse} twice for~$\bar{q} \coloneq q \circ \phi^{-1} \in C^{k}$ leaves us with the quadratic expression $\bar{q}(\bar{z}) = (q\circ\phi\inv)(\bar{z}) = \tfrac12 \bar{z}\tp \mathsf{Q}(\bar{z}) \bar{z}$ with pointwise symmetric~$\mathsf{Q} \in C^{k-2}$ satisfying $\mathsf{Q}(0) = \nabla^2 q(0) \succ 0$.
    The definiteness of~$\mathsf{Q}$ can be achieved by choosing~$U$ sufficiently small.
\end{proof}

Let~\eqref{eq:plant} be passive with storage function~$\ham$.
If \Cref{as:smoothness-hamc} is satisfied and~$\ham$ and~$\hamc$ satisfy the assumptions in \Cref{lem:morse},
then there exists an open set~$V \subseteq \RR^n$ with~$0\in V$ and a $C^k$-isomorphism~$\phi \colon \phi^{-1}(V) \to V, ~ z \mapsto \bar{z}$ with $\phi(0) = 0$ and
\begin{equation*}
    (\ham \circ \phi^{-1})(\bar{z}) = \frac12 \bar{z}\tp \bar{z}
    \und
    (\hamc \circ \phi^{-1})(\bar{z}) = \frac12 \bar{z}\tp \Hc(\bar{z}) \bar{z}\tp
\end{equation*} 
for pointwise symmetric positive definite $C^{k-2}$ function $\Hc \colon V \to \RR^{n,n}$ with $\Hc(0) = \Detac(0)$.
Diagonalizing $\Hc(\bar{z})$ yields 
\begin{equation}\label{eq:Hc-eigen-decomposition}
    \Hc(\bar{z}) = T(\bar{z}) D(\bar{z}) T(\bar{z})\tp,
\end{equation}
where the functions $D, T \colon V \to \RR^{n,n}$ collect the eigenvalues and eigenvectors, respectively, with $D$ pointwise diagonal and $T$ pointwise orthogonal.
To show that balanced coordinates exist, we will assume that $\theta \colon \bar{z} \mapsto T(\bar{z})\tp \bar{z}$ is bijective and of class $C^{k-2}$.
We are not aware of general results in this direction, and it appears that the assumption is not stated in~\cite[Theorem 4.3]{scherpen93-balancing} or the subsequent results~\cite{scherpen94-normalized,kawano18-structure}.
Under the assumption that the number of distinct eigenvalues is constant, the smoothness of the eigenvalues can be deduced from~\cite[Theorem 5.13a]{kato82-short}.
However, the reference only provides smoothness results for the eigenprojections rather than for the eigenvectors.
If the eigenvalues are pairwise distinct, then from~\cite[Equation II.4.1]{kato95-perturbation} it follows that the eigenvectors can be chosen as smooth as the eigenprojections.
To establish bijectivity of~$\theta$, the implicit function theorem may be utilized.

\begin{theorem}[balancing]\label{thm:balancing}
    In addition to \Cref{as:smoothness-hamc}, assume that~\eqref{eq:plant} is passive with storage function $\ham$, where~$\ham$ and~$\hamc$ satisfy the assumptions in \Cref{lem:morse} for $k\geq2$.
    Moreover, with~$T$ as in~\eqref{eq:Hc-eigen-decomposition}, assume that~$\theta \colon V \to \theta(V), ~ \bar{z} \mapsto T(\bar{z})\tp \bar{z}$ is a $C^{k-2}$-isomorphism.
    Then $\psi \coloneq \theta \circ \phi \colon \phi\inv(V) \to \theta(V), z \mapsto \tilde{z}$ is a $C^{k-2}$-isomorphism  with $\psi(0) = 0$ and 
    \begin{equation*}
        (\ham \circ \psi\inv)(\tilde{z}) = \frac12 \tilde{z}\tp \tilde{z}
        \und
        (\hamc \circ \psi\inv)(\tilde{z}) = \tilde{z}\tp \Sigma(\tilde{z}) \tilde{z},
    \end{equation*}
    where~$\Sigma \colon \theta(V) \to \RR^{n,n}$ is $C^{k-2}$ and pointwise diagonal with positive entries.
\end{theorem}
\begin{proof}
    Let $\tilde{z} \in \theta(V)$ be arbitrary.
    With the definition of~$\theta$ and the pointwise orthogonality of~$T$ we have $\theta\inv(\tilde{z})\tp \theta\inv(\tilde{z}) = \theta\inv(\tilde{z})\tp T(\tilde{z}) T(\tilde{z})\tp \theta\inv(\tilde{z}) = \tilde{z}\tp \tilde{z}$ and hence
    \begin{equation*}
        (\ham \circ \psi\inv)(\tilde{z})
        =
        (\ham \circ \phi\inv)(\theta\inv(\tilde{z}))
        =
        \tfrac12 \tilde{z}\tp \tilde{z}.
    \end{equation*}
    Moreover, setting $\Sigma \coloneq D \circ \theta\inv$ we obtain
    \begin{equation*}
        (\hamc \circ \psi\inv)(\tilde{z})
        =
        (\hamc \circ \phi\inv)(\theta\inv(\tilde{z}))
        =
        \frac12 \tilde{z}\tp \Sigma(\tilde{z}) \tilde{z}.
    \end{equation*}
    The other properties of $\phi$ and $\Sigma$ hold by definition.
\end{proof}

\Cref{thm:balancing} differs from~\cite[Theorem 3.14]{kawano18-structure} as follows.
First, in~\cite{kawano18-structure} the cost functional does not contain an output term.
In turn, the Hamilton--Jacobi--Bellman equation~\eqref{eq:control-hjb} simplifies, and choosing~$y=0$ in our result leads to results similar to~\cite{kawano18-structure}.
Second, in~\cite{kawano18-structure} the bijectivity assumption on~$\theta$ is not explicitly stated.
Third, we do not consider \portHamiltonian structure but general passive systems.
Lastly, \Cref{lem:morse} and \Cref{thm:balancing} do not assume the functions to be $C^\infty$.

To~\cite{sarkar23-structure}, the differences are the following.
In~\cite{sarkar23-structure}, three functions are balanced: a generalized controllability function~$\tilde{L}_c$, a generalized observability function~$\tilde{L}_o$, and the energy functional of a \portHamiltonian system, with $\tilde{L}_c$ being transformed to be in the form~$\tfrac12 z\tp z$.
In our setting,~$\hamc$ plays the role of the generalized observability function, we do not consider a generalized controllability function, and instead of \portHamiltonian systems we consider general passive systems.
Moreover, as mentioned before, we do not assume smoothness.

If the system~\eqref{eq:plant} is in coordinates as in \Cref{thm:balancing}, then we say that the system is \emph{semibalanced with respect to~$\ham$ and~$\hamc$}. 
Similar to~\cite{kawano18-structure,sarkar23-structure}, we can show that if the system is semibalanced with respect to~$\ham$ and~$\hamc$, then truncation preserves passivity.

\begin{theorem}[truncation]\label{thm:truncation}
    Let the assumptions from \Cref{thm:balancing} be satisfied.
    Further, assume that~\eqref{eq:plant} is semibalanced with respect to $\ham$ and $\hamc$, and that~\eqref{eq:plant} can be written in the form~\eqref{eq:r-definition}.
    If we partition $z = (\zk, \znk)$ with $\zk \in \RR^k$, $k \leq n$ arbitrary, and 
    \begin{equation*}
        r(z) = \vec{\rk(\zk, \znk) \\ \rnk(\zk, \znk)}
        \und
        g(z) = \vec{\gk(\zk, \znk) \\ \gnk(\zk,\znk)}
    \end{equation*}
    accordingly, then the reduced order model 
    \begin{equation}\label{eq:reduced-passive-system}
        \begin{aligned}
        \dot{z}_1 & = \rk(\zk, 0) + \gk(\zk, 0) u \\
        \yk & = \gk(\zk, 0)\tp \zk
        \end{aligned}
    \end{equation}
    is passive.
\end{theorem}
\begin{proof}
    In the semibalanced realization, we have $\ham(z) = \tfrac12 z\tp z$.
    Thus, plugging $v = z = (\zk, \znk)$ in~\eqref{eq:r-definition} with $\eta(z) = z$ gives
    \begin{equation*}
        0 \leq z\tp r(\eta(z)) 
        = \vec{\zk\tp & \znk\tp} \vec{\rk(\zk, \znk) \\ \rnk(\zk, \znk)} 
        = \zk\tp \rk(\zk, \znk) + \znk\tp \rnk(\zk, \znk).
    \end{equation*}
    Since $\znk = 0$ in~\eqref{eq:reduced-passive-system}, we obtain the claimed passivity.
\end{proof}

Similar to~\cite[Section III.A]{kawano18-structure}, where \portHamiltonian systems are considered, it can be shown that diffeomorphic state transformations preserve the structure~\eqref{eq:r-definition} if~$\eta$ is bijective.
To see that, let $\psi \colon \RR^n \to \RR^n$ be diffeomorphic and set $\tilde{\ham} \coloneq \ham \circ \psi\inv$.
Then $\nabla\tilde{\ham}(\tilde{z}) = D\psi\inv(\tilde{z})\tp \eta(\psi\inv(\tilde{z}))$ for all $\tilde{z} = \psi(z) \in \RR^n$, where $D\psi$ denotes the Jacobian of~$\psi$.
Similarly, $\dot{\tilde{z}} = D\psi(z) \dot{z}$.
By the inverse function rule
we have $D\psi\inv(\tilde{z}) = D\psi(z)\inv$.
Now, we can set $\tilde{r}(v) \coloneq D\psi(z) r(D\psi(z)\tp v) $ with $v \in \RR^n$ arbitrary and $z = \eta\inv(v)$, and combine the previous facts to establish the equivalence of 
\begin{equation*}
    \dot{z} = -r(\eta(z)) + g(z) u
    \und
    \dot{\tilde{z}} = - \tilde{r}(\nabla\tilde{\ham}(\tilde{z})) + \tilde{g}(\tilde{z}) u,
\end{equation*}
where $\tilde{g}(\tilde{z}) \coloneq D\psi(z) g(z)$ and $\tilde{z} = \psi(z)$.
Clearly, $v\tp \tilde{r}(v) \geq 0$ for all $v\in \RR^n$.
Hence, if the transformation~$\psi$ from \Cref{thm:balancing} is diffeomorphic and $\eta$ is bijective, then applying \Cref{thm:balancing,thm:truncation} consecutively yields a passivity preserving balanced truncation procedure.

The numerical computation of the balancing transformation from \Cref{thm:balancing} can be challenging.
Since the model order reduction scheme is not the main focus of this work, for an example of a practical realization of a related balanced truncation scheme the interested reader is referred to~\cite{sarkar23-structure,corbin24-scalable}.

\section{Numerical experiments}
\label{sec:numerical-experiments}

We are now ready to study the proposed controller~\eqref{eq:passive-controller} numerically.
We start by presenting our test examples, discuss the details of a practical realization, and finally study both the controller behavior and the proposed passivity.
In all experiments, we consider~$\Qcost = \Rcost = \eye$.

As a comparison for the controller performance, we use the Luenberger observer~\eqref{eq:luenberger} with the observer gain $L$ given by the extended Kalman filter~\cite{jazwinski70-stochastic}.
This leads to the controller 
\begin{equation}\label{eq:ekf-controller}
    \begin{aligned}
        \dzekf & = f(\zekf) - g(\zekf) g(\zekf)\tp \etac(\zekf) + K(\cdot, \zekf) (\uekf - h(\zekf)), \\
        \yekf & = g(\zekf) \etac(\zekf), \\
        \dot{\Pi} & = F(\zekf) \Pi + \Pi F(\zekf)\tp - \Pi H (\zekf)\tp \mathcal{R}\inv H(\zekf) \Pi + \mathcal{Q},
    \end{aligned}
\end{equation}
where $\zekf(0) = \zekf_0$ and $\Pi(0) = \Pi_0$ are the initial conditions, $F(\zekf) \coloneq D(f - gg\tp \etac)(\zekf)$, $H(\zekf) \coloneq Dh(\zekf)$, $K(t, \zekf(t)) \coloneq \Pi(t) H(\zekf(t))\tp \mathcal{R}\inv$ and $\mathcal{R} \in \RR^{m,m}$ and $\mathcal{Q} \in \RR^{n,n}$ are weighting matrices.
In all experiments we take $\mathcal{R} = I_m$, $\mathcal{Q} = I_n$ and $\zekf_0 = 0$ and $\Pi_0 = I_n$.

\subsection{Examples}
Since the numerical approximation of value functions is a challenging task for high dimensional spatial domains~\cite{kalise18-polynomial}, we focus on two-dimensional systems in our numerical examples.
The first example system models the motion of a nonlinear pendulum, and the corresponding dynamical system is passive.
The second system is a damped Van~der~Pol oscillator, where compared to the pendulum the nonlinearity has a stronger influence in the dynamics, making the numerical approximation of~$\hamc$ more challenging.

\begin{example}[nonlinear pendulum]\label{exmp:pendulum}
    The nonlinear system 
    \begin{equation*}
        \ddot \theta = - g \sin(\theta) - \lambda \dot\theta + u
    \end{equation*}
    describes the motion of a pendulum.
    Here, $\theta$ denotes the angular displacement of the pendulum, $g$ is the gravitational constant, $\lambda$ is a friction coefficient and $u$ is a control input.
    Using $z = (z_1, z_2) = (\theta, \dot\theta)$, we can express the system as 
    \begin{equation*}
        \vec{\dot{z_1} \\ \dot{z_2}}
        =
        \vec{
            z_2 \\
            - g \sin(z_1) - \lambda z_2 + u
        }
        =
        (J - R) \eta(z) + B u
    \end{equation*}
    with $J = \textvec{0 & 1 \\ -1 & 0}$, $R = \textvec{0 & 0 \\ 0 & \lambda}$, $B = \textvec{0 \\ 1}$ and $\eta = \nabla \ham$, where
    \begin{equation*}
        \ham(z) = g (1 - \cos(z_1)) + \tfrac12 z_2^2.
    \end{equation*}
    For the output we choose $y = B\tp \eta(z)$.
    In all numerical experiments we use $g = 9.81$ and $\lambda = 0.2$ and consider $z_0 = \textvec{\frac\pi4 & -1}\tp$.
\end{example}

\begin{example}[damped Van~der~Pol oscillator]\label{exmp:van_der_pol}
    We consider the forced and damped Van~der~Pol oscillator modeled by
    \begin{equation*}
        \ddot{x} = \mu (1 - x^2) \dot{x} - \lambda \dot{x} - x + u,
    \end{equation*}
    where $x$ denotes the position coordinate, $\mu$ and $\lambda$ are nonlinear and linear damping parameters, respectively, and $u$ is a forcing term.
    Using $z = (z_1, z_2) = (x, \dot{x})$, we can express the system as
    \begin{equation*}
        \vec{\dot{z_1} \\ \dot{z_2}}
        =
        \vec{
            z_2
            \\
            \mu (1 - z_1^2) z_2 - \lambda z_2 - z_1
        }
        +
        \vec{0 \\ 1} u
        = f(z) + B u.
    \end{equation*}
    In view of~\eqref{eq:ocp}, we consider the measurement $h(z) = \textvec{ 1 & 0} z$.
    In all numerical experiments we use $\mu=2$ and $\lambda=1.6$ and consider $z_0 = \textvec{1 & -0.5}\tp$.
\end{example}

\subsection{Practical realization}
\label{subsec:ocp-implementation}

For the practical realization of the controller~\eqref{eq:passive-controller}, we approximate the value function~$\hamc$ using the policy iteration~\cite{howard60-dynamic}.
The interested reader is referred to, e.g.,~\cite{kalise18-polynomial} for an overview on other available methods.
The policy iteration is based around the idea of taking an initial guess for the optimal policy $u^{(0)}$ and then iteratively solving the HJB equation~\eqref{eq:control-hjb} to obtain a new policy iterate via $u^{(k)} = - g(z)\tp \etac^{(k)}(z)$.
We focus on two-dimensional systems and utilize a Galerkin scheme on the rectangular space domain $\Omega = [\underline{x}, \overline{x}]\times[\underline{y},\overline{y}]$ based on scaled Legendre polynomials $\phi_i$, $\phi_j$, $i,j=1,\dots,d$, which form a $\lebesgue^2$-orthonormal basis of the polynomials of degree $\leq d$ on the intervals $[\underline{x},\overline{x}]$ and $[\underline{y},\overline{y}]$, respectively.
We make the ansatz 
\begin{equation}\label{eq:value-function-ansatz}
    \hamc^{(k+1)}(z) = \sum_{i,j=1}^{d} \alpha^{(k+1)}_{i,j} \psi_{i,j}(z),
\end{equation}
where $\psi_{i,j}(\vec{z_1 & z_2}) = \phi_i(z_1) \phi_j(z_2)$, $i,j = 1,\dots,d$ are the basis functions and $\alpha_{i,j}^{(k+1)}$ are basis coefficients.
After multiplying with a test function $\psi_{r,s}$ and integrating over the spatial domain, with $\etac^{(k+1)} = \nabla \hamc^{(k+1)}$ the HJB equation~\eqref{eq:control-hjb} becomes the linear system
\begin{equation}\label{eq:hjb-policy}
\begin{aligned}
    & \sum_{i,j=1}^{d} \alpha_{i,j}^{(k+1)}
    \left[ 
        \int_{\underline{x}}^{\overline{x}} 
        \int_{\underline{y}}^{\overline{y}}
            \nabla\psi_{i,j}\tp
            \big(f + B u^{(k)}\big)
            \psi_{r,s}
        \,\mathrm{d}z_2
        \,\mathrm{d}z_1
    \right]_{r,s=1,\dots,d}
    \\
    & =
    -
    \frac12
    \left[ 
        \int_{\underline{x}}^{\overline{x}} 
        \int_{\underline{y}}^{\overline{y}}
            \big(h\tp h + (u^{(k)})\tp u^{(k)}\big)
            \psi_{r,s}
        \,\mathrm{d}z_2
        \,\mathrm{d}z_1
    \right]_{r,s=1,\dots,d}\!\!,
\end{aligned}
\end{equation}
where we have suppressed the state dependencies for brevity.
For the initial policy, we linearize the plant~\eqref{eq:plant} around zero and choose $u^{(0)}(z) = - g(0)\tp \Pc z$, where $\Pc=\Pc\tp\succeq0$ is the stabilizing solution of~\eqref{eq:control-are}.
As a stopping criterion, we consider the quantities 
\begin{equation*}
    \delta^{(k+1)}_{\text{abs}} = \max_{z \in \Omegatest}{\norm{(u^{(k+1)} - u^{(k)})(z)}},~~
    \delta^{(k+1)}_{\text{rel}} = \frac{\delta^{(k+1)}_{\text{abs}}}{\max_{z \in \Omegatest}\norm{u^{(k)}(z)}},
\end{equation*}
where $\Omegatest \subseteq \Omega$ is the Cartesian product of $100$ equispaced samples in $x$ and $y$ direction.
The resulting method is summarized in \Cref{alg:policy-iteration}.

\begin{algorithm}
\caption{2D Policy iteration for~\eqref{eq:ocp} using Galerkin ansatz}\label{alg:policy-iteration}
\begin{algorithmic}[1]
    \Require 
        Initial control $u^{(0)}$, parameters $\tolabs$, $\tolrel$.
    \Ensure 
        Approximation of the value function $\hamc$, optimal control $u^*$
    \State
        Set $k=0$.
    \While{
        true
        }
        \State
            Compute $\hamc^{(k+1)}$ via~\eqref{eq:hjb-policy}, calculate $\etac^{(k+1)} = \nabla \hamc^{(k+1)}$.
        \State 
            Set $u^{(k+1)}(z) = - g(z)\tp \etac^{(k+1)}(z)$.
        \If{
            $\delta^{(k+1)}_{\text{abs}} \leq \tolabs$ 
            or
            $\delta^{(k+1)}_{\text{rel}} \leq \tolrel$
            }
            end while loop
        \EndIf
        \State
            Update $k \gets k+1$.
    \EndWhile
    \Return 
        $\hamc = \hamc^{(k)}$, $\etac = \etac^{(k)}$, $u^* = u^{(k)}$
\end{algorithmic}
\end{algorithm}

Unless stated otherwise, in all numerical experiments we initialize the controller with $\zc_0 = 0$.
For the stopping criterion, we choose $\tolabs = 10^{-14}$ and $\tolrel = 10^{-10}$.
In all of our experiments, the policy iteration converged in under 10 steps.

\subsection{Controller behavior}

To study the behavior of the controller, we consider the systems from \Cref{exmp:pendulum,exmp:van_der_pol} as plants, where for the maximal degree of the one dimensional ansatz polynomials we chose $d=10$ and $d=15$, respectively.
Since passivity of the closed loop system is not guaranteed, in all experiments the trajectories of the coupled dynamics were computed using the implicit midpoint method with $500$ time steps.

In \Cref{fig:pendulum_value_function}, the value function computed by the policy iteration is shown.
Since \Cref{exmp:pendulum} is well approximated by a linear system, the value function appears to be quadratic.
\Cref{fig:pendulum_coupled_trajectory} depicts the coupled trajectory of~\eqref{eq:plant} and \eqref{eq:passive-controller} for \Cref{exmp:pendulum}.
To study the effectiveness of choosing $B$ as the observer gain in~\eqref{eq:luenberger}, in \Cref{fig:pendulum_coupled_trajectory_ekf} the coupled trajectory is shown when the controller~\eqref{eq:passive-controller} is replaced by~\eqref{eq:ekf-controller}.
As shown in \Cref{fig:pendulum_state_decay}, in both cases the behavior of the controlled plant clearly differs from the behavior of the uncontrolled plant, steering the plant to zero significantly faster.
We observe that for this example, switching to the controller~\eqref{eq:ekf-controller} accelerates convergence to zero and increases the magnitude of the controller state.

In \Cref{fig:van_der_pol_value_function,fig:van_der_pol_coupled_trajectory_ekf}, the same experiments are shown for \Cref{exmp:van_der_pol}.
Here, the value function is clearly not quadratic, reflecting the strong influence of the nonlinearity in \Cref{exmp:van_der_pol}.
As before, compared to the uncontrolled plant dynamics, using the controller~\eqref{eq:passive-controller} accelerates the convergence to zero of the plant trajectory.
Similar to the results for \Cref{exmp:pendulum}, we observe that the switching to the controller~\eqref{eq:ekf-controller} leads to larger values of the controller state.
Moreover, in contrast to the previous experiment, the convergence to zero of both the plant and the controller trajectory is slowed down by the switch to~\eqref{eq:ekf-controller}.


\begin{figure}
    \centering
	\begin{subfigure}[t]{.45\textwidth}
        \centering
        \scalebox{0.48}{\input{results/pendulum_value_function.pgf}}
        \caption{
            Value function of~\eqref{eq:ocp}.
        }
        \label{fig:pendulum_value_function}
    \end{subfigure}
	\hfill
	\begin{subfigure}[t]{.45\textwidth}
		\centering
        \scalebox{0.48}{\input{results/pendulum_state_decay.pgf}}
        \caption{
            Decay of the plant trajectory for the uncontrolled case and the coupled system with~\eqref{eq:passive-controller} and~\eqref{eq:ekf-controller}.
        }
        \label{fig:pendulum_state_decay}
	\end{subfigure}
    \\
    \begin{subfigure}[t]{.45\textwidth}
		\centering
        \scalebox{0.48}{\input{results/pendulum_coupled_trajectory.pgf}}
        \caption{
            Dynamics of the coupling of the plant~\eqref{eq:plant} and the controller~\eqref{eq:passive-controller}.
        }
        \label{fig:pendulum_coupled_trajectory}
	\end{subfigure}
    \hfill
    \begin{subfigure}[t]{.45\textwidth}
		\centering
        \scalebox{0.48}{\input{results/pendulum_coupled_trajectory_ekf.pgf}}
        \caption{
            Dynamics of the coupling of the plant~\eqref{eq:plant} and the controller~\eqref{eq:ekf-controller}.
        }
        \label{fig:pendulum_coupled_trajectory_ekf}
	\end{subfigure}
    \caption{
        Results for the nonlinear pendulum (\Cref{exmp:pendulum}).
    }
\end{figure}


\begin{figure}
    \centering
	\begin{subfigure}[t]{.45\textwidth}
        \centering
        \scalebox{0.48}{\input{results/van_der_pol_value_function.pgf}}
        \caption{
            Value function of~\eqref{eq:ocp}.
        }
        \label{fig:van_der_pol_value_function}
    \end{subfigure}
	\hfill
    \begin{subfigure}[t]{.45\textwidth}
		\centering
        \scalebox{0.48}{\input{results/van_der_pol_state_decay.pgf}}
        \caption{
            Decay of the plant trajectory of for the uncontrolled case and the coupled system with~\eqref{eq:passive-controller} and~\eqref{eq:ekf-controller}.
        }
        \label{fig:van_der_pol_state_decay}
	\end{subfigure}
    \\
    \begin{subfigure}[t]{.45\textwidth}
		\centering
        \scalebox{0.48}{\input{results/van_der_pol_coupled_trajectory.pgf}}
        \caption{
            Dynamics of the coupling of the plant~\eqref{eq:plant} and the controller~\eqref{eq:passive-controller}.
        }
        \label{fig:van_der_pol_coupled_trajectory}
	\end{subfigure}
    \hfill
    \begin{subfigure}[t]{.45\textwidth}
		\centering
        \scalebox{0.48}{\input{results/van_der_pol_coupled_trajectory_ekf.pgf}}
        \caption{
            Dynamics of the coupling of the plant~\eqref{eq:plant} and the controller~\eqref{eq:ekf-controller}.
        }
        \label{fig:van_der_pol_coupled_trajectory_ekf}
	\end{subfigure}
    \caption{
        Results for the damped Van~der~Pol oscillator (\Cref{exmp:pendulum}).
    }
\end{figure}

\subsection{Passivity verification}

It is well-known that classical time discretization methods can interfere with the passivity structure of dynamical systems~\cite{hairer06-geometric}.
For instance, the implicit midpoint rule preserves quadratic invariants of dynamical systems, but may not preserve other invariants.
Hence, to numerically verify the passivity of the controller~\eqref{eq:passive-controller}, a passivity preserving time discretization scheme is needed.
Different structure-preserving methods are available, and we refer to~\cite{giesselmann25-energy} for an overview of recent developments.
Here, we use a time discretization method based on the discrete gradient framework~\cite{harten83-upstream,itoh88-hamiltonian,gonzalez96-time}.
See also~\cite[Section V.5]{hairer06-geometric} and~\cite[Page 1024]{mclachlan99-geometric} for historical remarks.

Consider a system of the form~\eqref{eq:plant} and assume that it is passive with the continuously differentiable storage function $\ham$. 
Taking the time derivative in~\eqref{eq:energybalance}, we obtain the power balance
\begin{equation}\label{eq:powerbalance}
    \ddt \ham(z) = \eta(z)\tp f(z) + y\tp u \leq y\tp u.
\end{equation}
Our goal is to state a time discretization scheme that qualitatively preserves the power balance and dissipation inequality~\eqref{eq:powerbalance}, see~\eqref{eq:time-discrete-powerbalance} for the precise notion.
Our approach is similar to the one taken in~\cite{mclachlan99-geometric}.
Therein, the authors assume that the dynamics may be expressed in the linear gradient form~$f(z) = L(z) \eta(z)$, where~$L \colon \RR^n \to \RR^{n,n}$ satisfies~$L(z) + L(z)\tp \preceq 0$ for all~$z\in\RR^n$.
Since the controller~\eqref{eq:passive-controller} may not be written in this form in general, see the discussion at the end of \Cref{sec:passive-controller}, their method is not applicable for our setting.
Nevertheless, similar ideas may be used.

A continuous function $\etabar \colon \RR^n \times \RR^n \to \RR$ is called \emph{discrete gradient of $\ham$} if $\etabar$ satisfies the mean value property
\begin{equation}\label{eq:discrete-gradient-property-eta}
    \ham(z_2) - \ham(z_1) = \etabar(z_1, z_2)\tp(z_1 - z_2)
\end{equation}
and the consistency property
\begin{equation*}
    \etabar(z_1,z_1) = \eta(z_1)
\end{equation*} 
for all~$z_1, z_2 \in \RR^n$.
If $n=1$, then the unique discrete gradient is the difference quotient.
In higher dimensions, the discrete gradient is not unique, since only the component along~$z_2 - z_1$ is restricted.
As noted by~\cite{mclachlan99-geometric}, this limits the ability of discrete gradients to approximate point values of continuous gradients to second order.
High order extensions are discussed in, e.g.,~\cite{eidnes22-order}, see also the references therein.
In what follows, we will use the \emph{Gonzalez discrete gradient}~\cite{gonzalez96-time}, which reads
\begin{equation}\label{eq:def-etabar}
    \etabar(z_1,z_2) =
    \begin{cases}
        \eta\!\left( \frac{z_1 + z_2}{2} \right) 
        + 
        \frac{
                \ham(z_2) 
                - \ham(z_1) 
                - \eta\left( \frac{z_1 + z_2}{2} \right)\tp (z_2 - z_1) 
            }{
                \norm{z_2 - z_1}^2
            } 
            (z_2 - z_1) 
        \\ \hfill \text{if $z_1 \neq z_2$,}
        \\
        \eta(z_1) \hfill \text{else.}
    \end{cases}
\end{equation}

In addition to discrete gradients, we employ the representation~\eqref{eq:r-definition}, which motivates the following discretization scheme.

\begin{theorem}[passivity preserving time discretization]\label{thm:discrete-gradient-method}
    Let time points $0 = t_{1} < \dots < t_{m} = T$ be given and define $\gbar(z_1, z_2) \coloneq g\! \left( \frac{z_1 + z_2}{2} \right)$ for $z_1, z_2 \in \RR^n$ as well as $\ubar_{i} \coloneq \frac{u(t_{i}) + u(t_{i+1})}{2}$ for $i = 1, \dots, m-1$.
    Further, assume~$r$ satisfying~\eqref{eq:r-definition} is given.
    Then a solution of the time-discrete system 
    \begin{equation}\label{eq:time-discrete-system}
            z_{i+1} 
            =
        z_{i} - (t_{i+1} - t_{i}) r\big(\etabar(z_{i}, z_{i+1})\big)
        + (t_{i+1} - t_{i}) \gbar(z_{i}, z_{i+1}) \ubar_{i}, 
    \end{equation}
    $i = 1, \dots, m-1$,
    satisfies the time-discrete power balance and dissipation inequality
    \begin{equation}\label{eq:time-discrete-powerbalance}
            \frac{\ham(z_{i+1}) - \ham(z_{i})}{t_{i+1} - t_{i}} 
            = 
            -
            \etabar(z_{i}, z_{i+1})\tp 
            r\big(\etabar(z_{i}, z_{i+1})\big)
            + \ybar_i\tp \ubar_i
            \leq
            \ybar_i\tp \ubar_i,
    \end{equation}
    where the discrete output~$\ybar_i$ is defined as $\ybar_i \coloneq \gbar(z_{i}, z_{i+1})\tp \etabar(z_{i}, z_{i+1})$, $i=0,\dots,m-1$.
\end{theorem}
\begin{proof}
    Rearranging~\eqref{eq:time-discrete-system} yields 
    \begin{equation*}
        \frac{z_{i+1} - z_{i}}{t_{i+1} - t_{i}} 
        = 
        -r\big(\etabar(z_{i}, z_{i+1})\big)
        + 
        \gbar(z_{i}, z_{i+1}) \ubar_i.
    \end{equation*}
    Using~\eqref{eq:discrete-gradient-property-eta} and inserting the equation above, we obtain 
    \begin{align*}
        \frac{\ham(z_{i+1}) - \ham(z_{i})}{t_{i+1} - t_{i}} 
        & =
        \etabar(z_{i}, z_{i+1})\tp \frac{z_{i+1} - z_{i}}{t_{i+1} - t_{i}} 
        \\
        & =
        - \etabar(z_{i},z_{i+1})\tp r\big(\etabar(z_{i}, z_{i+1})\big) 
        + \etabar(z_{i},z_{i+1})\tp \gbar(z_{i}, z_{i+1}) \ubar_i.
    \end{align*}
    Using that $\etabar(z_1,z_2)\tp r\big(\etabar(z_1, z_2)\big) \geq 0$ for all $z_1, z_2 \in \RR^n$ by~\eqref{eq:r-definition} finishes the proof.
\end{proof}

Before we verify the passivity of~\eqref{eq:passive-controller} using the time discretization scheme discussed in \Cref{thm:discrete-gradient-method}, we investigate the performance of the scheme using \Cref{exmp:pendulum} as a test scenario.
In all numerical experiments, we considered a constant step size $\Delta t \coloneq t_{i+1} - t_i$.

In \Cref{fig:pendulum_convergence}, the convergence of the scheme for \Cref{exmp:pendulum} with control $u(t) = \sin(t)$ is shown.
The reference solution was computed with the implicit midpoint rule with step size $\Delta t_{\text{ref}} = 2^{-3} \cdot \Delta t_{0}$, where $\Delta t_{0} = 10^{-3}$ is the finest shown step size in \Cref{fig:pendulum_convergence}.
For the implicit midpoint rule, we approximated $u(\frac{t_{i} + t_{i+1}}{2}) \approx \frac{u(t_{i}) + u(t_{i+1})}{2}$.
We observe that the proposed discrete gradient scheme appears to be second order accurate.
In \Cref{fig:pendulum_powerbalance}, the power balance~\eqref{eq:time-discrete-powerbalance} is verified numerically by visualizing the relative error
\begin{equation*}
    \frac{
        \left|
            \frac{\ham(z_{i+1}) - \ham(z_{i})}{t_{i+1} - t_{i}}
            + \etabar(z_{i}, z_{i+1})\tp r(\etabar(z_{i}, z_{i+1})) 
            - \ybar_i\tp \ubar_i
        \right|
        }{  
            \max_{j = 1,\dots,m}
            \left|
                \frac{\ham(z_{j+1}) - \ham(z_{j})}{t_{j+1} - t_{j}}
            \right|
        }.
\end{equation*}
We use $m=500$ and observe that the relative error is close to machine precision throughout the time horizon.

\begin{figure}
	\begin{subfigure}[t]{.49\textwidth}
		\centering
		\scalebox{0.48}{\input{results/pendulum_convergence.pgf}}
		\caption{
            Convergence of the solution with $u(t) = \sin(t)$.
		}
		\label{fig:pendulum_convergence}
	\end{subfigure}
	\hfill
	\begin{subfigure}[t]{.49\textwidth}
		\centering
		\scalebox{0.48}{\input{results/pendulum_powerbalance.pgf}}
		\caption{
            Relative error in the power balance~\eqref{eq:time-discrete-powerbalance} with $u(t) = \sin(t)$.
        }
		\label{fig:pendulum_powerbalance}
	\end{subfigure}
    \caption{
        Discrete gradient method applied to the nonlinear pendulum (\Cref{exmp:pendulum}).
    }
\end{figure}

To numerically verify the passivity of the controller~\eqref{eq:passive-controller}, we consider \Cref{exmp:pendulum,exmp:van_der_pol}.
As before, for the maximal degree of the ansatz polynomials in~\eqref{eq:value-function-ansatz}, we chose $d=10$ and $d=15$, respectively.
For the computation of $\etac\inv(v)$ for $r$ as in~\eqref{eq:r-definition}, we used a Newton iteration for $g(\zc) \coloneq \etac(\zc) - v$.
We used the initial guess $\zc = 0$ and performed $10$ steps of the Newton iteration, which yielded $\norm{g(\zc)}$ in machine precision for the final iterate.
In \Cref{fig:pendulum_discrete_gradient,fig:van_der_pol_discrete_gradient}, the behavior of $\hamc(\zc(t))$ for \eqref{eq:passive-controller} with control input $\uc(t) = 0$ and initial state $\zc(0) = \textvec{1 & 1}\tp$ and the relative error in the power balance~\eqref{eq:time-discrete-powerbalance} is shown for \Cref{exmp:pendulum,exmp:van_der_pol}.
In both cases, we consider the time horizon $[0,10]$ and use the time stepping scheme~\eqref{eq:time-discrete-system} with $m=500$.
As expected, we observe that $\hamc(\zc(t))$ is non-increasing over the length of the time horizon for both examples, and the relative error in the power balance is close to machine precision for the entirety of the time horizon.

\begin{figure}
	\begin{subfigure}[t]{.49\textwidth}
		\centering
		\scalebox{0.48}{\input{results/pendulum_controller_hamiltonian_zero_control.pgf}}
		\caption{
            Behavior of $\hamc(\zc(t))$.
        }
		\label{fig:pendulum_controller_hamiltonian_zero_control}
	\end{subfigure}
    \hfill
    \begin{subfigure}[t]{.49\textwidth}
		\centering
		\scalebox{0.48}{\input{results/pendulum_controller_powerbalance_zero_control.pgf}}
		\caption{
            Relative error in the power balance~\eqref{eq:time-discrete-powerbalance}.
        }
		\label{fig:pendulum_controller_powerbalance_zero_control}
	\end{subfigure}
    \caption{
        Discrete gradient method applied to the controller~\eqref{eq:passive-controller} for the nonlinear pendulum (\Cref{exmp:pendulum}) with $\uc(t) = 0$ and $\zc(0) = \textvec{1 & 1}\tp$.
    }
    \label{fig:pendulum_discrete_gradient}
\end{figure}

\begin{figure}
	\begin{subfigure}[t]{.49\textwidth}
		\centering
		\scalebox{0.48}{\input{results/van_der_pol_controller_hamiltonian_zero_control.pgf}}
		\caption{
            Behavior of $\hamc(\zc(t))$
        }
		\label{fig:van_der_pol_controller_hamiltonian_zero_control}
	\end{subfigure}
    \hfill
    \begin{subfigure}[t]{.49\textwidth}
		\centering
		\scalebox{0.48}{\input{results/van_der_pol_controller_powerbalance_zero_control.pgf}}
		\caption{
            Relative error in the power balance~\eqref{eq:time-discrete-powerbalance}.
        }
		\label{fig:van_der_pol_controller_powerbalance_zero_control}
	\end{subfigure}
    \caption{
        Discrete gradient method applied to the controller~\eqref{eq:passive-controller} for the damped Van~der~Pol oscillator (\Cref{exmp:van_der_pol}) with $\uc(t) = 0$ and $\zc(0) = \textvec{1 & 1}\tp$.
    }
    \label{fig:van_der_pol_discrete_gradient}
\end{figure}

\section{Conclusion}\label{sec:conclusion}
In this paper, we developed a passive controller for nonlinear systems with affine input structure.
Our work extends the results of~\cite{breiten22-error} to the nonlinear setting.
Furthermore, we discussed applications of our results for nonlinear balancing based model order reduction.
Finally, using a novel passivity preserving time discretization scheme, we numerically verified the passivity of the controller and investigated its effect on multiple nonlinear examples.
Here, we observed that the controller has a similar performance as a feedback controller using the extended Kalman filter for the observer gain.

Possible topics for future research include the study of passive controllers with robustness guarantees, extending the results of~\cite{breiten23-structure}, and investigating applications of the model order reduction scheme from \Cref{sec:passive-controller} utilizing results of~\cite{sarkar23-structure}.

\paragraph{Acknowledgments}
The authors thank the Deutsche Forschungsgemeinschaft for their support within the subproject B03 in the Sonderforschungsbereich/Transregio 154 “Mathematical Modelling, Simulation and Optimization using the Example of Gas Networks” (Project 239904186).
The authors would also like to thank P.\,Schulze and J.\,Schröder for their helpful comments and discussions, 
as well as the anonymous reviewers for their helpful feedback and suggestions for improvement.

\bibliographystyle{plain}        
\bibliography{phcon}           

\begin{thebibliography}{10}

\bibitem{altmann25-structure}
R.~Altmann, A.~Karsai, and P.~Schulze.
\newblock Structure-preserving discretization and model reduction for energy-based models.
\newblock {\em arXiv preprint {2507.21552}}, 2025.

\bibitem{antoulas05-approximation}
A.~C. Antoulas.
\newblock {\em Approximation of Large-Scale Dynamical Systems}.
\newblock Society for Industrial and Applied Mathematics, 2005.

\bibitem{breiten23-structure}
T.~Breiten and A.~Karsai.
\newblock Structure-preserving {$\mathcal{H}_{\infty}$} control for port-{Hamiltonian} systems.
\newblock {\em Systems {\&} Control Letters}, 174:105493, 2023.

\bibitem{breiten22-error}
T.~Breiten, R.~Morandin, and P.~Schulze.
\newblock Error bounds for port-{Hamiltonian} model and controller reduction based on system balancing.
\newblock {\em Computers {\&} Mathematics with Applications}, 116:100--115, 2022.

\bibitem{breiten25-structure}
T.~Breiten and P.~Schulze.
\newblock Structure-preserving linear quadratic gaussian balanced truncation for port-{Hamiltonian} descriptor systems.
\newblock {\em Linear Algebra and its Applications}, 704:146--191, 2025.

\bibitem{breiten22-passivity}
T.~Breiten and B.~Unger.
\newblock Passivity preserving model reduction via spectral factorization.
\newblock {\em Automatica}, 142:110368, 2022.

\bibitem{buchfink23-symplectic}
P.~Buchfink, S.~Glas, and B.~Haasdonk.
\newblock Symplectic model reduction of {Hamiltonian} systems on nonlinear manifolds and approximation with weakly symplectic autoencoder.
\newblock {\em SIAM Journal on Scientific Computing}, 2023.

\bibitem{buchfink24-model}
P.~Buchfink, S.~Glas, B.~Haasdonk, and B.~Unger.
\newblock Model reduction on manifolds: A differential geometric framework.
\newblock {\em Physica {D}: Nonlinear Phenomena}, 468:134299, 2024.

\bibitem{camlibel23-port}
M.~K. Camlibel and A.~J. Van~der Schaft.
\newblock Port-{Hamiltonian} systems theory and monotonicity.
\newblock {\em {SIAM} Journal on Control and Optimization}, 2023.

\bibitem{coddington55-theory}
Earl~A. Coddington and Norman Levinson.
\newblock {\em Theory of Ordinary Differential Equations}.
\newblock International series in pure and applied mathematics. McGraw-Hill, New York, NY, USA, 1955.

\bibitem{corbin24-scalable}
N.~A. Corbin, A.~Sarkar, J.~M.~A. Scherpen, and B.~Kramer.
\newblock Scalable computation of input-normal/output-diagonal balanced realization for control-affine polynomial systems.
\newblock {\em arXiv preprint {2410.22435}}, 2024.

\bibitem{crandall83-viscosity}
M.~G. Crandall and P.-L. Lions.
\newblock Viscosity solutions of {Hamilton}-{Jacobi} equations.
\newblock {\em Transactions of the American Mathematical Society}, 277:1--42, 1983.

\bibitem{eidnes22-order}
S{\o}lve Eidnes.
\newblock Order theory for discrete gradient methods.
\newblock {\em BIT Numerical Mathematics}, 62:1207--1255, 2022.

\bibitem{fleming06-controlled}
W.~H. Fleming and H.~M. Soner.
\newblock {\em Controlled Markov Processes and Viscosity Solutions}.
\newblock Springer, New York, NY, USA, 2006.

\bibitem{freeman96-robust}
R.~A. Freeman and P.~Kokotovi\'{c}.
\newblock {\em Robust Nonlinear Control Design}.
\newblock Birkhäuser, Boston, MA, USA, 1996.

\bibitem{gernandt25-port}
H.~Gernandt and M.~Schaller.
\newblock Port-{Hamiltonian} structures in infinite-dimensional optimal control: {P}rimal--{D}ual gradient method and control-by-interconnection.
\newblock {\em Systems {\&} Control Letters}, 197:106030, 2025.

\bibitem{giesselmann25-energy}
J.~Giesselmann, A.~Karsai, and T.~Tscherpel.
\newblock Energy-consistent {Petrov}--{Galerkin} time discretization of port-{Hamiltonian} systems.
\newblock {\em {SMAI} Journal of computational mathematics}, 11:335--367, 2025.

\bibitem{gonzalez96-time}
O.~Gonzalez.
\newblock Time integration and discrete {Hamiltonian} systems.
\newblock {\em Journal of Nonlinear Science}, 6:449--467, 1996.

\bibitem{gugercin12-structure}
S.~Gugercin, R.~Polyuga, C.~Beattie, and A.~{van der Schaft}.
\newblock Structure-preserving tangential interpolation for model reduction of {port-Hamiltonian} systems.
\newblock {\em Automatica}, 48(9):1963 -- 1974, 2012.

\bibitem{haddad08-nonlinear}
W.~M. Haddad and V.~Chellaboina.
\newblock {\em Nonlinear Dynamical Systems and Control: A {Lyapunov}-Based Approach}.
\newblock Princeton University Press, Princeton, NJ, USA, 2008.

\bibitem{hairer06-geometric}
E.~Hairer, G.~Wanner, and C.~Lubich.
\newblock {\em Geometric Numerical Integration}.
\newblock Springer, Berlin, G., 2006.

\bibitem{harten83-upstream}
A.~Harten, P.~D. Lax, and B.~van Leer.
\newblock On upstream differencing and {Godunov}-type schemes for hyperbolic conservation laws.
\newblock {\em SIAM Review}, 1983.

\bibitem{hill80-dissipative}
D.~J. Hill and P.~J. Moylan.
\newblock Dissipative dynamical systems: Basic input-output and state properties.
\newblock {\em Journal of the Franklin Institute}, 309:327--357, 1980.

\bibitem{howard60-dynamic}
R.~A. Howard.
\newblock Dynamic programming and markov processes.
\newblock 2:39--47, 1960.

\bibitem{isidori99-nonlinear}
A.~Isidori.
\newblock {\em Nonlinear Control Systems {II}}.
\newblock Springer, London, England, UK, 1999.

\bibitem{itoh88-hamiltonian}
T.~Itoh and K.~Abe.
\newblock Hamiltonian-conserving discrete canonical equations based on variational difference quotients.
\newblock {\em Journal of Computational Physics}, 76:85--102, 1988.

\bibitem{jazwinski70-stochastic}
Andrew~H Jazwinski.
\newblock {\em Stochastic processes and filtering theory}.
\newblock Academic Press, 1970.

\bibitem{jin21-inverse}
W.~Jin, Dana Kuli\'{c}, S.~Mou, and S.~Hirche.
\newblock Inverse optimal control from incomplete trajectory observations.
\newblock {\em International Journal of Robotics Research}, 40:848--865, 2021.

\bibitem{johnson79-state}
C.~Johnson.
\newblock State-variable design methods may produce unstable feedback controllers.
\newblock {\em International Journal of Control}, 29(4):607--619, 1979.

\bibitem{kalise18-polynomial}
D.~Kalise and K.~Kunisch.
\newblock Polynomial approximation of high-dimensional {Hamilton}--{Jacobi}--{Bellman} equations and applications to feedback control of semilinear parabolic {PDE}s.
\newblock {\em SIAM Journal on Scientific Computing}, 2018.

\bibitem{kalman64-when}
R.~E. Kalman.
\newblock When is a linear control system optimal?
\newblock {\em Journal of Basic Engineering}, 86:51--60, 1964.

\bibitem{karsai25-passivity}
A.~Karsai, T.~Breiten, J.~Ramme, and P.~Schulze.
\newblock Passivity encoding representations of nonlinear systems.
\newblock {\em IEEE Transactions on Automatic Control}, 70:7660--7666, 2025.

\bibitem{kato82-short}
T.~Kato.
\newblock {\em A Short Introduction to Perturbation Theory for Linear Operators}.
\newblock Springer, New York, NY, USA, 1982.

\bibitem{kato95-perturbation}
T.~Kato.
\newblock {\em Perturbation Theory for Linear Operators}.
\newblock Springer, Berlin, Germany, 2nd edition, 1995.

\bibitem{kawano18-structure}
Y.~Kawano and J.~M.~A. Scherpen.
\newblock Structure preserving truncation of nonlinear port {Hamiltonian} systems.
\newblock {\em IEEE Transactions on Automatic Control}, 63:4286--4293, 2018.

\bibitem{lang85-differential}
S.~Lang.
\newblock {\em Differential Manifolds}.
\newblock Springer, New York, NY, USA, 1985.

\bibitem{luenberger64-observing}
D.~G. Luenberger.
\newblock Observing the state of a linear system.
\newblock {\em IEEE Transactions on Military Electronics}, 8:74--80, 1964.

\bibitem{mclachlan99-geometric}
R.~I. McLachlan, G.~R.~W. Quispel, and N.~Robidoux.
\newblock Geometric integration using discrete gradients.
\newblock {\em Philosophical Transactions of the Royal Society of London. Series {A}: Mathematical, Physical and Engineering Sciences}, 357:1021--1045, 1999.

\bibitem{mehrmann23-control}
V.~Mehrmann and B.~Unger.
\newblock Control of port-{Hamiltonian} differential-algebraic systems and applications.
\newblock {\em Acta Numerica}, 32:395--515, 2023.

\bibitem{milnor16-morse}
J.~Milnor.
\newblock {\em {Morse} Theory}.
\newblock Princeton University Press, Princeton, NJ, USA, 2016.

\bibitem{moylan74-implications}
P.~Moylan.
\newblock Implications of passivity in a class of nonlinear systems.
\newblock {\em IEEE Transactions on Automatic Control}, 19:373--381, 1974.

\bibitem{moylan73-nonlinear}
P.~Moylan and B.~Anderson.
\newblock Nonlinear regulator theory and an inverse optimal control problem.
\newblock {\em IEEE Transactions on Automatic Control}, 18:460--465, 1973.

\bibitem{ortega89-adaptive}
R.~Ortega and M.~W. Spong.
\newblock Adaptive motion control of rigid robots: A tutorial.
\newblock {\em Automatica}, 25:877--888, 1989.

\bibitem{ortega02-interconnection}
R.~Ortega, A.~{Van der Schaft}, B.~Maschke, and G.~Escobar.
\newblock Interconnection and damping assignment passivity-based control of port-controlled {H}amiltonian systems.
\newblock {\em Automatica}, 38(4):585--596, apr 2002.

\bibitem{sarkar23-structure}
A.~Sarkar and J.~M.~A. Scherpen.
\newblock Structure-preserving generalized balanced truncation for nonlinear port-{Hamiltonian} systems.
\newblock {\em Systems {\&} Control Letters}, 174:105501, 2023.

\bibitem{scherpen93-balancing}
J.~M.~A. Scherpen.
\newblock Balancing for nonlinear systems.
\newblock {\em Systems {\&} Control Letters}, 21:143--153, 1993.

\bibitem{scherpen94-normalized}
J.~M.~A. Scherpen and A.~J. Van Der~Schaft.
\newblock Normalized coprime factorizations and balancing for unstable nonlinear systems.
\newblock {\em International Journal of Control}, 1994.

\bibitem{schwerdtner23-sobmor}
P.~Schwerdtner and M.~Voigt.
\newblock {SOBMOR}: Structured optimization-based model order reduction.
\newblock {\em SIAM Journal on Scientific Computing}, 2023.

\bibitem{sepulchre97-constructive}
R.~Sepulchre, M.~Janković, and P.~V. Kokotović.
\newblock {\em Constructive Nonlinear Control}.
\newblock Springer, London, England, UK, 1997.

\bibitem{vanderschaft17-l2gain}
A.~{Van der Schaft}.
\newblock {\em $L^2$-Gain and Passivity Techniques in Nonlinear Control}, volume~2.
\newblock Springer International Publishing, 2017.

\bibitem{vanderschaft14-port}
A.~{Van der Schaft} and D.~Jeltsema.
\newblock {Port-Hamiltonian} systems theory: An introductory overview.
\newblock {\em Foundations and Trends in Systems and Control}, 1(2-3):173--378, 2014.

\bibitem{wang03-generalized}
Y.~Wang, C.~Li, and D.~Cheng.
\newblock Generalized {Hamiltonian} realization of time-invariant nonlinear systems.
\newblock {\em Automatica}, 39:1437--1443, 2003.

\bibitem{willems72-dissipative2}
J.~Willems.
\newblock Dissipative dynamical systems part {II}: Linear systems with qua\-drat\-ic supply rates.
\newblock {\em Archive for Rational Mechanics and Analysis}, 45(5):352--393, 1972.

\bibitem{willems72-dissipative1}
Jan Willems.
\newblock Dissipative dynamical systems part {I}: general theory.
\newblock {\em Archive for Rational Mechanics and Analysis}, 45(5):321--351, 1972.

\bibitem{zanon21-constrained}
M.~Zanon and A.~Bemporad.
\newblock Constrained controller and observer design by inverse optimality.
\newblock {\em IEEE Transactions on Automatic Control}, 67:5432--5439, 2021.

\end{thebibliography}

\end{document}